%% file: main.tex
\documentclass[10pt]{article} 
\include{0Prae}
 
\title{Approximation in the extended functional tensor train format}
\usepackage{authblk}
\author[1]{Christoph Strössner}
\author[1]{Bonan Sun}
\author[1]{Daniel Kressner}
\affil[1]{\small \'Ecole Polytechnique F\'ed\'erale de Lausanne (EPFL), Institute of Mathematics, CH-1015 Lausanne, Switzerland.  \url{christoph.stroessner@epfl.ch}, \url{bonan.sun@epfl.ch}, \url{daniel.kressner@epfl.ch}}
\date{February 23, 2024}

\begin{document}
\maketitle

\begin{abstract} 
This work proposes the extended functional tensor train (EFTT) format for compressing and working with multivariate functions on tensor product domains. Our compression algorithm combines tensorized Chebyshev interpolation with a low-rank approximation algorithm that is entirely based on function evaluations. Compared to existing methods based on the functional tensor train format, \rerevision{the adaptivity of} our approach often \rerevision{results in reducing} the required storage, sometimes considerably,  while achieving the same accuracy. In particular, we reduce the number of function evaluations required to achieve a prescribed accuracy by up to over $96\%$ compared to the algorithm from [Gorodetsky, Karaman and Marzouk, \emph{Comput. Methods Appl. Mech. Eng.}, 347 (2019)] . 
\end{abstract}

\input{1Introduction}
\input{2Setting}

\input{3Algorithm}
\input{4NumericalComparison}
\input{5Conclusion}

\subsection*{Acknowledgments}
The authors would like to thank Sergey Dolgov and \rerevision{both reviewers for their very helpful comments.}
The authors declare that they have no conflict of interest.
{\footnotesize

\input{main.bbl}
}

\input{6Apendix}

\end{document}

%% file: 0Prae.tex
\usepackage{fullpage}
\usepackage[T1]{fontenc} 
\usepackage[latin1]{inputenc}
\usepackage[english]{babel}
\usepackage{selinput}
\SelectInputMappings{adieresis={ä},germandbls={ß}}

 \usepackage{hyperref}

\usepackage{amsmath} 
\usepackage{amssymb} 
\usepackage{amsthm} 
\usepackage{enumerate}

\newtheorem{lemma}{Lemma}

\newtheorem{rmrk}{Remark}

\usepackage{graphicx} 
\usepackage{epstopdf}
\usepackage{booktabs}
\usepackage{multirow}
\usepackage{subfig}
\usepackage[framed,numbered]{matlab-prettifier}
\usepackage{algorithm} 
\usepackage[noend]{algpseudocode}

\graphicspath{{./Figures/}}

\usepackage{multirow}
\usepackage{fancyvrb}
\usepackage{hyperref}

\newcommand{\R}{\mathbb{R}}

\newcommand{\norm}[1]{\lVert #1 \rVert}
\newcommand{\abs}[1]{\left\lvert #1 \right\rvert}

\newcommand{\eps}{\varepsilon}

\DeclareMathOperator{\argmax}{argmax}

\newcommand{\uargmax}[1]{\underset{#1}{\argmax}\;}

\newcommand{\umax}[1]{\underset{#1}{\max}\;}

\usepackage{ mathrsfs }

\newcommand{\T}{\mathcal{T}}
\usepackage[numbers,sort]{natbib}
\usepackage{svg}

\newcommand{\revision}[1]{\textcolor{black}{ #1}}
\newcommand{\rerevision}[1]{\textcolor{black}{ #1}}

%% file: 1Introduction.tex
\section{Introduction}

Multivariate functions become notoriously difficult to approximate and work with in high dimensions. The storage and computational effort of standard approximation techniques increases exponentially with the dimension. Several techniques have been developed to potentially mitigate this so-called curse of dimensionality, including sparse grids~\cite{Bungartz04} and low-rank tensor approximation techniques~\cite{Hackbusch12}. On a functional level, low-rank tensor formats correspond to separation of variables. In particular, if $f:[-1,1]^d \to \mathbb R$ can be written as a sum of separable functions,
\begin{equation} \label{eq:cpfun}
  f(x_1,\dots,x_d) = \sum_{\alpha = 1}^R g_{\alpha}^{(1)}(x_1) g_{\alpha}^{(2)}(x_2)\cdots g_{\alpha}^{(d)}(x_d),  
\end{equation}
then the task of approximating $f$ is replaced by approximating the $R \cdot d$ univariate functions $g_{\alpha}^{(\ell)}: [-1,1] \to \mathbb R$; see, e.g.,~\cite{Beylkin02}. It is often beneficial, especially when $R$ is very large, to impose additional structure on~\eqref{eq:cpfun}, leading to the so-called functional tensor train (FTT)~\cite{Bigoni16,Gorodetsky19} and Tucker~\cite{Hashemi17,Stroessner21} formats. 
In practice, functions are rarely given directly in such a functional low-rank format, but they can sometimes be well approximated by it.
This is exploited in quite a few applications, including the solution of time-dependent partial differential equations (PDEs)~\cite{Dektor19},  uncertainty quantification~\cite{Konakli16}, sensitivity analysis~\cite{Ripoll19}, optimal control~\cite{Gorodetsky18b} and quantum dynamic simulations~\cite{Soley21}.

The FTT format~\cite{Bigoni16,Gorodetsky19} takes the form
\begin{equation}
    \label{eq:FunTTformat}
    f(x_1,\dots,x_d) = \sum_{\alpha_1 = 1}^{R_1} \dots, \sum_{\alpha_{d-1}=1}^{R_{d-1}} g^{(1)}_{1,\alpha_1}(x_1) g^{(2)}_{\alpha_1,\alpha_2}(x_2) \cdots g^{(d-1)}_{\alpha_{d-2},\alpha_{d-1}}(x_{d-1}) g^{(d)}_{\alpha_{d-1},1}(x_d),
\end{equation}
for univariate functions $g^{(\ell)}_{\alpha_{\ell - 1},\alpha_\ell}:[-1,1] \to \mathbb R$ with $\alpha_{\ell - 1} =1,\dots,R_{\ell -1}$, $ \alpha_\ell = 1,\dots, R_{\ell}$,  $\ell = 1,\dots, d$. The summation ranges $(R_1,\dots,R_{d-1}) \in \mathbb N^{d-1}$ are called TT (representation) ranks and we formally set $R_0 = R_d = 1$ to simplify notation.  From~\eqref{eq:FunTTformat}, a fully discrete approximations is obtained via discretizing each $g^{(\ell)}_{\alpha_{\ell - 1},\alpha_\ell}$ by, e.g., a truncated series expansion, which incurs additional truncation error~\cite{Gorodetsky19}.

It is important to emphasize that not every function admits a good approximation in the FTT format~\eqref{eq:FunTTformat} with moderate TT ranks, 
especially for larger $d$. Generally speaking, smoothness~\cite{Hackbusch07} and a restricted (e.g., nearest neighbor) interaction between variables are helpful.
Upper bounds for the TT ranks needed to attain a certain accuracy are derived in~\cite{Bigoni16,Griebel19} for functions in Sobolev spaces. 
The required TT ranks can change significantly when variables are transformed~\cite{Trefethen17c} or the order of variables is permuted~\cite{Dektor22}. 
For functions in periodic and mixed Sobolev spaces approximation rates can be found in~\cite{Schneider14} and~\cite{Griebel22}, respectively.
Much better rates can be obtained for functions with special structures such as compositional functions~\cite{Bachmayr21} and quantities associated with certain parametric PDEs~\cite{Bachmayr17}.
For a more abstract analysis of functional low-rank approximation we refer to~\cite{Ali21,Ali20b,Ali20a}.

To compute functional low-rank approximations in practice, one often starts from multivariate polynomial interpolation~\cite{Mason80,Trefethen17b}. 
When using a tensorized grid of interpolation nodes, the number of function evaluations and the memory needed to store the coefficients of the interpolant grow exponentially with $d$. This growth can be alleviated by utilizing a low-rank approximation of the coefficient tensor, which corresponds to a functional low-rank approximation~\cite{Hackbusch12,Stroessner21}. Such an approach is used for bivariate~\cite{Townsend13} and trivariate~\cite{Hashemi17,Stroessner21} functions by Chebfun~\cite{Driscoll14} and for general multivariate functions in~\cite{Bigoni16}. 
The impact of a TT approximation of the coefficient tenor on the overall interpolation error has been studied for tensorized Lagrange interpolation with respect to the $L^2$-norm in~\cite{Bigoni16}.
In~\cite{Stroessner21}, the impact of arbitrary coefficient tensor approximations is analyzed with respect to the uniform norm in the context of tensorized Chebyshev interpolation of trivariate functions.
A slightly different approach to obtain FTT approximations is proposed in~\cite{Gorodetsky19}, which
applies a continuous version of the so-called TT-cross algorithm~\cite{Oseledets10} directly to the function.
All univariate functions occurring in this process are discretized using parameterizations such as basis expansions.
In principle, the format allows every univariate function to be parameterized differently.
In the special case that univariate interpolation with the same basis functions is used to discretize all univariate functions in the same mode, this approach is equivalent to a TT approximation of the coefficient tensor for tensorized interpolation.

The FTT format~\eqref{eq:FunTTformat} requires the storage of $R_{\ell-1} \cdot R_{\ell}$ univariate functions for each $\ell = 1,\dots,d$.
These functions (or their parametrizations) will often be linearly dependent.
In this case, we can compress the format further by storing these functions in terms of linear combinations of $ r_{\ell} \ll R_{\ell-1} \cdot R_{\ell}$ basis functions.
The resulting approximation can be seen as a variant of the hierarchical Tucker format proposed in~\cite{Hackbusch09}.
When the approximation is constructed from a low-rank approximation of the coefficient tensor, we can obtain such a compressed format by following the ideas of Khoromskij~\cite{Khoromskij06}, who 
suggested to use a two-level Tucker format, which is obtained by first computing a Tucker approximation~\cite{Tucker66} before approximating the core tensor  in $\R^{r_1 \times \dots \times r_d}$ further. 
Computing a TT approximation of the core tensor leads to a so-called extended TT format~\cite{Eigel20,Dolgov13,Schneider14}, from which we can derive a functional low-rank approximation in the extended functional tensor train (EFTT) format. This EFTT approximation corresponds to a compressed FTT approximation.

In this work, we propose a novel algorithm to efficiently compute functional low-rank approximations using the EFTT format.
Our algorithm is based on first obtaining suitable factor matrices for a Tucker approximation of the coefficient tensor. 
The columns of these matrices are determined by applying a combination of adaptive cross approximation~\cite{Bebendorf00} with randomized pivoting to matricizations of the coefficient tensor.
We then apply oblique projections based on discrete empirical interpolation~\cite{Chaturantabut10} as in~\cite{Sorensen16,Stroessner21} to implicitly construct a suitable core tensor for the Tucker approximation.
In a second step, we compute a TT approximation of this core tensor using the greedy restricted cross interpolation algorithm~\cite{Savostyanov14}, which is a rank adaptive variant of the TT-cross algorithm~\cite{Oseledets10} requiring asymptotically fever evaluations than the TT-DMRG-cross algorithm~\cite{Savostyanov11}.
The combination of the Tucker and TT approximation yields the desired functional low-rank approximation.
The main advantage of our approach compared to a direct TT approximation of the coefficient tensor is that the approximated core tensor is potentially much smaller than the coefficient tensor.
This reduces the storage complexity of the approximation and it significantly reduces the number of function evaluations needed in the computation of the TT approximation, which we demonstrate in our numerical experiments. \rerevision{Other approaches that can potentially achieve this goal have been discussed in~\cite{Grelier18,HaberstichPhD,Bertrand22,trunschke2023weighted}.}

\revision{Our novel algorithm is potentially useful in all applications that require to store and work with a multivariate functions. This includes uncertainty quantification, computational chemistry and physics; we refer to the surveys~\cite{Grasedyck13,Khoromskaia18,Khoromskij18} for concrete examples.}

\revision{The TT format considered in this work corresponds to a line graph tensor network topology. Other tensor network topologies have been considered in the literature for the purpose of function approximation. In~\cite{Ballani13}, a variant of the TT-cross algorithm for tree tensor networks, corresponding to the hierarchical Tucker format, has been developed. In~\cite{Haberstich21}, strategies for choosing a good tree tensor network have been proposed. It would be interesting to compare and possible combine the developments from~\cite{Haberstich21} with this work.} 

The remainder of this paper is structured as follows. 
In Section~\ref{sec:problemsetting}, we recall how functional low-rank approximations can be obtained by combining tensorized Chebyshev interpolation and low-rank approximations of the coefficient tensor.
Our novel algorithm to compute approximations in the EFTT format is presented in Section~\ref{sec:Algorithm}.
In our numerical experiments in Section~\ref{sec:NumericalExperiments}, we demonstrate the advantages of our novel algorithm and compare it to a direct TT approximation of the coefficient tensor as in~\cite{Bigoni16} and to the FTT algorithm in~\cite{Gorodetsky19}.

%% file: 2Setting.tex
\section{Functional low-rank approximation via tensorized interpolation}\label{sec:problemsetting}

\subsection{Tensorized Chebyshev interpolation} \label{sec:tencheb}

Given a function $f:[-1,1]^d\to \R$, let us consider a polynomial approximation of degree $(n_1,\dots,n_d)$ taking the form
\begin{equation}\label{eq:ChebyshevInterpolantForm}
   f(x_1,\dots,x_d) \approx \tilde{f}(x_1,\dots,x_d) = \sum_{i_1=0}^{n_1} \dots \sum_{i_d=0}^{n_d} \mathcal{A}_{i_1,\dots,i_d} T_{i_1}(x_1)  \cdots T_{i_d}(x_d),
\end{equation}
where $\mathcal{A}\in \R^{(n_1+1) \times \dots \times (n_d+1)}$ is the coefficient tensor and $T_k(x) = \cos(k\cos^{-1}(x))$ denotes the $k$-th Chebyshev polynomial. To be consistent with standard notation~\cite{NIST}, we index entries of the coefficient tensor starting from $0$. Matrices and tensors not related to Chebyshev approximation are indexed in the usual way, starting from $1$.

It is common to construct the approximation~\eqref{eq:ChebyshevInterpolantForm} by interpolating $f$ on a tensorized grid of Chebyshev points~\cite{Trefethen13}, yielding a tensor $\mathcal{T}\in \R^{(n_1+1)\times\cdots\times (n_d+1)}$ containing the function values:
\begin{equation}
\label{eq:EvaluationTensor}
	\mathcal{T}_{i_1,\dots,i_d}=f\big(x_{i_1}^{(1)},\dots,x_{i_d}^{(d)}\big),\  x_{k}^{(\ell)}=\cos(\pi k / n_\ell),\ k=0,\dots,n_\ell,\  \ell=1,\dots, d.
\end{equation}
The coefficient tensor $\mathcal{A}$ in~\eqref{eq:ChebyshevInterpolantForm} is constructed by
applying discrete cosine transformations~\cite{Gentleman72}.
This can be written as multiplication of $\mathcal{T}$ with matrices $\mathcal{F}^{(\ell)}$ encoding these transformations:
\begin{equation*}
	\mathcal{A}=\mathcal{T}\times_{1} F^{(1)}\times_2 F^{(2)} \times _3\dots\times_d F^{(d)},
\end{equation*}
where $\times_\ell$ denotes the mode-$k$ multiplication~\cite{Kolda09} defined as
\[\mathcal{A}=\mathcal{T}\times_{\ell} F^{(\ell)} \iff \mathcal{A}_{i_1,\dots,i_d}=\sum_{j_\ell=0}^{n_\ell}\mathcal{T}_{i_1,\dots,i_{\ell-1},j_\ell,i_{\ell+1},\dots,i_d} F_{i_\ell,j_\ell}^{(\ell)} \text{ for all } i_\ell = 0,\dots,n_\ell,\ \ell = 1,\dots,d\]
and $F^{(\ell)} \in \R^{(n_\ell + 1)\times (n_\ell + 1)}$ is \revision{derived from}~\cite[Equation~6.28]{Mason02}:
\begin{equation*}
F^{({ \ell })} = \frac{2}{n_{{ \ell }}}
\begin{pmatrix}
\frac{1}{4} T_0(x_0^{({ \ell })}) & \frac{1}{2} T_0(x_1^{({ \ell })}) & \frac{1}{2} T_0(x_2^{({ \ell })}) & \dots & \frac{1}{4} T_{0}(x_{n_{{ \ell }} }^{({ \ell })}) \\
\frac{1}{2} T_1(x_0^{({ \ell })}) & T_1(x_1^{({ \ell })}) &  T_1(x_2^{({ \ell })}) & \dots & \frac{1}{2} T_{1}(x_{n_{{ \ell }} }^{({ \ell })}) \\
\frac{1}{2} T_2(x_0^{({ \ell })}) &  T_2(x_1^{({ \ell })}) &  T_2(x_2^{({ \ell })}) & \dots & \frac{1}{2} T_{2}(x_{n_{{ \ell }} }^{({ \ell })})
\\ \vdots & \vdots & \vdots & \ddots & \vdots \\
\frac{1}{4} T_{n_{{ \ell }}}(x_0^{({ \ell })}) & \frac{1}{2} T_{n_{{ \ell }}}(x_1^{({ \ell })}) & \frac{1}{2} T_{n_{{ \ell }}}(x_2^{({ \ell })}) & \dots & \frac{1}{4} T_{n_{{ \ell }}}(x_{n_{{ \ell }} }^{({ \ell })})
\end{pmatrix}.
\end{equation*}
For analytic functions $f$ it is shown in~\cite[Lemma 7.3.3]{Sauter11} that the error of the approximation~\eqref{eq:ChebyshevInterpolantForm} in the uniform norm decays exponentially with respect to $\min\{n_1,\ldots,n_d\}$.

\subsection{Functional tensor train (FTT) approximation}

The cost of constructing and storing the tensors $\mathcal{T}$ and $\mathcal{A}$ from Section~\ref{sec:tencheb}
 grows exponentially with $d$. To mitigate this growth, we replace $\mathcal{T}$ by a low-rank approximation $\hat{\mathcal{T}}$. The following lemma is a direct generalization of~\cite[Lemma~2.2]{Stroessner21} bounding the error introduced by such an approximation.

\begin{lemma} \label{lemma:approxlowrank}
Let $\tilde{f}$ be defined as in~\eqref{eq:ChebyshevInterpolantForm}. For  $\hat{\mathcal{T}} \in \mathbb R^{(n_1+1) \times \dots \times (n_d+1)}$, consider the polynomial 
\begin{equation}
\label{eq:HatF}
    \hat f(x_1,\dots,x_n) = \sum_{i_1=0}^{n_1} \dots \sum_{i_d=0}^{n_d} \hat{\mathcal{A}}_{i_1,\dots,i_d} T_{i_1}(x_1)  \cdots T_{i_d}(x_d),
\end{equation}
with $\hat{\mathcal{A}} = \hat{\mathcal{T}} \times_{1} F^{(1)}\times_2 \dots\times_d F^{(d)}$. Then
\[ 
\norm{f-\hat{f}}_{\infty} \leq \norm{f-\tilde{f}}_{\infty} + \prod_{\ell = 1}^{d} \big(\frac{2}{\pi}\log(n_\ell )+1\big) \norm{\T - \hat{\T}}_{\infty},
\]
where $\norm{\cdot}_\infty$ denotes the uniform norm for functions and the maximum norm for tensors. 
\end{lemma}

In the following, we approximate the tensor $\mathcal T$ by a tensor $\hat{\mathcal{T}}$ in TT format~\cite{Oseledets11}, the discrete analogue of~\eqref{eq:FunTTformat}. The entries of $\hat{\mathcal{T}}$ take the form
\begin{equation}
    \label{eq:TTformat}
    \hat{\mathcal{T}}_{i_1,\dots,i_d} = \sum_{\alpha_1 = 1}^{R_1} \dots \sum_{\alpha_{d-1}=1}^{R_{d-1}}  \mathcal G^{(1)}_{1,i_1,\alpha_1} \mathcal G^{(2)}_{\alpha_1,i_1,\alpha_2} \cdots \mathcal G^{(d-1)}_{\alpha_{d-2},i_{d-1},\alpha_{d-1}} \mathcal G^{(d)}_{\alpha_{d-1},i_d,1},
\end{equation}
with the so-called TT cores $\mathcal G^{(\ell)} \in \R^{R_{\ell - 1}\times (n_\ell+1) \times R_\ell}$ for $\ell = 1,\ldots,d$. If the maximal TT rank $R = \max R_\ell$ remains modest, this representation requires a storage of size $\mathcal{O}((d-2)nR^2 + 2nR)$ with $n = \max n_\ell$. This results in linear (instead of exponential) growth with respect $d$ under the (strong) assumption that $R$ remains constant as $d$ increases. 

Given $\hat{\mathcal{T}}$ in the TT format~\eqref{eq:TTformat}, the tensor $\hat{{\mathcal{A}}}$ defined in Lemma~\ref{lemma:approxlowrank} can also be expressed in TT format with the TT cores $\mathcal{A}^{(\ell)}:=\mathcal{G}^{(\ell)} \times_2 F^{(\ell)}$. In turn,
the function $\hat{f}$ defined in\eqref{eq:HatF} can be expressed in the FTT format~\eqref{eq:FunTTformat}, with the corresponding univariate functions given by 
\begin{equation*}
    g_{\alpha_{\ell - 1},\alpha_{\ell}}^{(\ell)}(x) = \sum_{j = 0}^{n_{\ell}}\mathcal A^{(\ell)}_{\alpha_{\ell - 1},j,\alpha_{\ell}} T_j(x) = \sum_{j=0}^{n_{\ell}} \sum_{k=0}^{n_{\ell}} F^{(\ell)}_{j,k} \mathcal G^{(\ell)}_{\alpha_{\ell - 1},k,\alpha_{\ell}} T_j(x).
\end{equation*}

\subsection{Extended functional tensor train (EFTT) approximation}
\begin{figure}[!ht]
    \centering
    \includegraphics[scale=1]{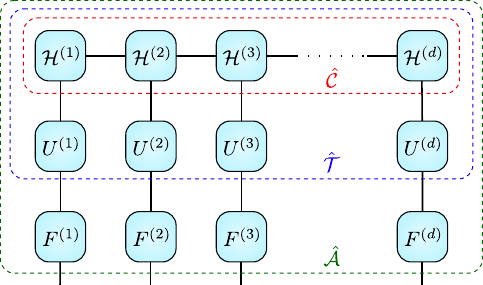}
    \caption{Tensor network representation~\cite{Orus14} of the coefficient tensor $\hat{\mathcal{A}}$ in~\eqref{eq:CoeffTensorExtendedFunTT} corresponding to EFTT format~\eqref{eq:ExtendedFunctionalTT}. The colored boxes mark the subtensors corresponding to the approximation of the evaluation tensor $\hat{\mathcal{T}}$ as in~\eqref{eq:HatF} and the approximation of the core tensor $\hat{\mathcal{C}}$ as in~\eqref{eq:TTofCore}.}
    \label{fig:TensorNetwork}
\end{figure}

As discussed in the introduction, we can potentially compress FTT approximations further by essentially replacing the TT core $\mathcal{G}^{(\ell)}$  in~\eqref{eq:TTformat} by $\mathcal{H}^{(\ell)}\times_2 U^{(\ell)}$, where $\mathcal{H}^{(\ell)} \in \R^{R_{\ell-1} \times r_\ell \times \R_{\ell}}$ and $U^{(\ell)} \in \R^{(n_\ell+1) \times r_\ell}$ with $r_\ell \ll R_{\ell-1} \cdot R_{\ell}$.
To compute such a compressed approximation without explicitly forming and then compressing a FTT approximation, we propose to compute an approximation of the evaluation tensor~\eqref{eq:EvaluationTensor} in the extended TT format~\cite{Schneider14}.
In this work, we obtain such an approximation by first approximating the evaluation tensor $\mathcal{T}$ from~\eqref{eq:EvaluationTensor} in the Tucker format~\cite{Tucker66}:
\begin{equation}\label{eq:Tucker}
    \mathcal{T} \approx \mathcal{C} \times_1 U^{(1)}\times_2 \dots \times_d U^{(d)},
\end{equation}
where $\mathcal{C}\in \R^{r_1\times \dots \times r_d}$ is called the core tensor, $U^{(\ell)} \in \R^{(n_\ell+1)\times r_{\ell}}$ are called factor matrices. The multilinear rank of $\mathcal{T}$ is given by $r_1,\ldots,r_{d}$ (more precisely, the smallest values for these integers such that $\mathcal{T}$ admits~\eqref{eq:Tucker}).
In a second step, the core tensor $\mathcal C$ is approximated in TT format:
\begin{equation}\label{eq:TTofCore}
    \mathcal{C}_{i_1,\dots,i_d} \approx  \hat{\mathcal{C}}_{i_1,\dots,i_d}  = \sum_{\alpha_1 = 1}^{R_1} \cdots \sum_{\alpha_{d-1}=1}^{R_{d-1}}  \mathcal H^{(1)}_{1,i_1,\alpha_1} \mathcal H^{(2)}_{\alpha_1,i_1,\alpha_2} \cdots \mathcal H^{(d-1)}_{\alpha_{d-2},i_{d-1},\alpha_{d-1}} \mathcal H^{(d)}_{\alpha_{d-1},i_d,1}
\end{equation}
with the TT cores $\mathcal{H}^{(\ell)}\in \R^{R_{\ell-1}\times r_{\ell}\times R_{\ell}}$.
Inserted into~\eqref{eq:Tucker}, this yields an approximation in extended TT format:
\begin{equation}\label{eq:ExtendedTTformat}
    \mathcal{T}_{i_1,\dots,i_d} \approx
    \sum_{j_1 = 1}^{r_1} \dots \sum_{j_d = 1}^{r_d} \sum_{\alpha_1 = 1}^{R_1} \dots \sum_{\alpha_{d-1}}^{R_{d-1}} \mathcal H^{(1)}_{1,j_1,\alpha_1}  \cdots  \mathcal H^{(d)}_{\alpha_{d-1},j_d,1} U^{(1)}_{i_1,j_1}\cdots U^{(d)}_{i_d,j_d}.
\end{equation}
This only requires $\mathcal{O}(drR^2+dnr)$ storage, where $R = \max R_\ell$, $r = \max r_\ell$, $n= \max n_\ell$, which compares favorably with the $\mathcal{O}(dnR^2)$ storage needed by the TT approximation~\eqref{eq:TTformat} alone, especially when when $r \ll n$.
From the approximation~\eqref{eq:ExtendedTTformat} of $\mathcal{T}$, we obtain the coefficient tensor approximation
\begin{equation}
    \label{eq:CoeffTensorExtendedFunTT}
\hat{\mathcal{A}}_{i_1,\dots,i_d} = \sum_{k_1=0}^{n_1} \dots \sum_{k_d=0}^{n_d} \sum_{j_1 = 1}^{r_1} \dots \sum_{j_d = 1}^{r_d} \sum_{\alpha_1 = 1}^{R_1} \dots \sum_{\alpha_{d-1}}^{R_{d-1}} \mathcal H^{(1)}_{1,j_1,\alpha_1}  \cdots  \mathcal H^{(d)}_{\alpha_{d-1},j_d,1}  U^{(1)}_{i_1,j_1}\cdots U^{(d)}_{i_d,j_d} F_{k_1,i_1}^{(1)} \cdots F_{k_d,i_d}^{(d)}.
\end{equation}
See Figure~\ref{fig:TensorNetwork} for a visualization of this approximation. 
Given $\hat{\mathcal{A}}$, we obtain an approximation of $f$ in the EFTT format:
\begin{equation}
    \label{eq:ExtendedFunctionalTT}
    f(x_1,\dots,x_d) \approx \sum_{j_1 = 1}^{r_1} \dots \sum_{j_d = 1}^{r_d} \sum_{\alpha_1 = 1}^{R_1} \dots \sum_{\alpha_{d-1}}^{R_{d-1}} \mathcal H^{(1)}_{1,j_1,\alpha_1}  \cdots  \mathcal H^{(d)}_{\alpha_{d-1},j_d,1} u^{(1)}_{j_1}(x_1)\cdots u^{(d)}_{j_d}(x_d),
\end{equation}
with the univariate functions $u^{(\ell)}_j: [-1,1] \to \R$ defined as $u^{(\ell)}_j(x) := \sum_{i = 0}^{n_\ell} \sum_{k=0}^{n_\ell} F_{k,i}^{(\ell)}  U_{i,j}^{(\ell)} T_k(x)$ for $1\leq \ell \leq d$, $1 \leq j \leq r_\ell$. 
\revision{To compute point evaluations of~\eqref{eq:ExtendedFunctionalTT}, we precompute $(F^{(\ell)}U^{(\ell)}) \in \R^{(n_\ell+1)\times r_\ell}$ using the discrete cosine transform for $\ell = 1,\dots,d$. Then, the evaluation at $(x_1,\dots,x_d) \in [-1,1]^d$ can be computed} by first contracting $F^{(\ell)}U^{(\ell)}$  with the vectors $(T_0(x_\ell),\dots,T_{n_\ell}(x_\ell))$ for $\ell = 1,\dots,d$ and by then contracting the resulting vectors with the TT cores $\mathcal{H}^{(\ell)}$~\cite{Orus14}.
This requires $\mathcal{O}(drn+drR^2)$ operations \revision{for each point evaluation}, where $R = \max R_\ell$, $r = \max r_\ell$, $n= \max n_\ell$. \revision{Additionally, $\mathcal{O}(dn \log(n))$ operations are required once to precompute the discrete cosine transforms.}

%% file: 3Algorithm.tex
\section{Approximation algorithm} \label{sec:Algorithm}

In the following, we develop a novel algorithm for computing approximations in the EFTT format~\eqref{eq:ExtendedFunctionalTT}. Our algorithm obtains the factor matrices $U^{(1)},\ldots,U^{(d)}$ in the Tucker approximation~\eqref{eq:Tucker} from fibers of the evaluation tensor $\mathcal T$
via a variant of column subset selection~\cite{Deshpande10,Cortinovis19}. Following~\cite{Stroessner21}, the core tensor $\mathcal C$ is obtained as a subtensor of $\mathcal T$ by applying discrete empirical interpolation. Thus, there is no need to form $\mathcal C$ explicitly and we apply a variant of the TT-cross algorithm~\cite{Oseledets10} to compute the TT cores $\mathcal H^{(1)},\dots,\mathcal{H}^{(d)}$ for the TT approximation~\eqref{eq:TTofCore} from only some entries of $\mathcal{C}$.

\subsection{Factor matrices}
\paragraph{Fiber approximations.} When fixing all but the $\ell$th index of $\mathcal T$, one obtains a vector called a mode-$\ell$ fiber of $\mathcal X$. The mode-$\ell$ matricization of $\mathcal{T}$
is a matrix $\mathcal{T}^{\{\ell\}} \in \R^{(n_\ell+1) \times ((n_1+1) \cdots (n_{\ell-1}+1) (n_{\ell+1}+1) \cdots (n_d+1))}$ that collects all mode-$\ell$ fibers~\cite{Kolda09}.
Applying adaptive cross approximation (ACA)~\cite{Bebendorf00}, a popular low-rank approximation technique that accesses a matrix via its entries,
to $\mathcal{T}^{\{\ell\}}$ yields index sets $\hat I_\ell ,\hat J_\ell$
that determine an approximation of the form 
\begin{equation}
\label{eq:ACAform}
    \mathcal{T}^{\{\ell\}} \approx \mathcal{T}^{\{\ell\}}(:,\hat J_\ell) (\mathcal{T}^{\{\ell\}}(\hat I_\ell,\hat J_\ell))^{-1} \mathcal{T}^{\{\ell\}}(\hat I_\ell,:).
\end{equation}
Note that the matrix $\mathcal{T}^{\{\ell\}}(:,\hat J_\ell)$ contains mode-$\ell$ fibers of $\mathcal{T}$.
If the approximation error of~\eqref{eq:ACAform} is small, we choose $\hat U^{(\ell)} = \mathcal{T}^{\{\ell\}}(:,\hat J_\ell)$ as an approximate basis of mode-$\ell$ fibers.

ACA with full pivoting determines the indices in $\hat I_\ell,\hat J_\ell$ subsequently by choosing the entry of largest absolute value in the current approximation error. This
corresponds to greedily maximizing the volume of $\mathcal{T}^{\{\ell\}}(\hat I_\ell,\hat J_\ell)$~\cite{Goreinov97} and comes with theoretical guarantees~\cite{Cortinovis20b}.
However, it is impossible to apply full pivoting in our context because it requires the evaluation of all entries. Even partial 
pivoting~\cite{Bebendorf03} is impractical as it requires the evaluation of a (very long) row in $\mathcal{T}^{\{\ell\}}$ in every step. Inspired by the success of randomization in numerical linear algebra~\cite{Martinsson20,Woodruff14}, we propose to determine the next indices based on sampling.
This leads to Algorithm~\ref{alg:ACA}, which samples a fixed number $s$ of random entries and picks the entry of largest absolute value. 
These random entries are also used for stopping the algorithm and determine the approximation rank $r_\ell$ adaptively. Note that the update in line~\ref{line:ACArecursiveEvaluation} is only performed implicitly and the subtraction is carried out each time an entry of $A$ is evaluated. 
Overall, Algorithm~\ref{alg:ACA} applied to $\mathcal{T}^{\{\ell\}}$ determines the index sets $\hat I_\ell,\hat J_\ell$ using only $\mathcal{O}(r_\ell^3+sr_\ell^2)$ entries of $\mathcal T$.
Afterwards, we explicitly compute $\hat U^{(\ell)} = \mathcal{T}^{\{\ell\}}(:,\hat J_\ell)$ by evaluating $n_\ell r_\ell$ entries of $\mathcal{T}$.
We perform the described procedure for every $\ell = 1,\dots,d$ to determine  $\hat U^{(1)},\dots,\hat U^{(d)}$; see also lines~\ref{line:ComputeJhat} and~\ref{line:ComputeUhat} of Algorithm~\ref{alg:Tucker} below.

\begin{algorithm}[!ht]
\caption{Adaptive cross approximation with randomized pivoting}\label{alg:ACA} 
\begin{algorithmic}[1]
\State \textbf{Input:} procedure to evaluate entries of $M\in\R^{n\times m}$, tolerance $\eps$, number of samples $s$
\State \textbf{Output:} index sets $I$ and $J$ such that $M \approx M(:,J)M(I,J)^{-1}M(I,:)$
\State $I = \emptyset$, $J = \emptyset$
\State $A = M$
\While {true}
\State Construct $S \subset \{1,\dots,n\}\times\{1,\dots,m\}$ by uniformly sampling $s$ index pairs $(i,j)$. \label{line:tupleSampling}
\If{$\umax{(i,j)\in S} \abs{A_{i,j}} \leq \eps$} return $I,J$ {\bf end} \EndIf
\State $(i^*,j^*) = \uargmax{(i,j)\in S} \abs{A_{i,j}}$
\State $I = I \cup \{i^*\}$, $J = J \cup \{j^*\}$
\State $A = M - M(:,J)M(I,J)^{-1}M(I,:)$ \label{line:ACArecursiveEvaluation}
\EndWhile
\end{algorithmic}
\end{algorithm}

\paragraph{Tucker approximation}
To arrive at a Tucker approximation~\eqref{eq:Tucker} we need to project the fibers of $\mathcal T$  onto the spans of $\hat U^{(1)},\dots,\hat U^{(d)}$. 
\revision{Let $\hat U^{(\ell)} = Q^{(\ell)} R^{(\ell)}$ denote the (economic) QR decomposition, where $Q^{(\ell)} \in \R^{(n_\ell+1) \times r_\ell}$ has orthonormal columns.}
\revision{The orthogonal projections used in the higher order singular value decomposition (HOSVD)~\cite{Lathauwer00}, i.e. multiplying $\mathcal T$ with $Q^{(\ell)}(Q^{(\ell)})^T$ in each mode}, would require to \revision{explicitly construct the} whole tensor $\mathcal T$. 
To avoid this, we follow the ideas of~\cite{Stroessner21,Sorensen16}. Given \revision{a} set $I_\ell$ containing $r_\ell$ indices selected from $\{0,\dots,n_\ell\}$, we let $\Phi_{I_\ell}$ contain the corresponding columns of the \revision{$(n_\ell+1) \times (n_\ell+1)$} identity matrix, such that 
$\Phi_{I_\ell}^T Q^{(\ell)}  = Q^{(\ell)}{(I_\ell,:)}$. If this matrix is invertible then
\[
 Q^{(\ell)} (\Phi_{I_\ell}^T Q^{(\ell)} )^{-1} \Phi_{I_\ell}^T
\]
defines an oblique projection onto the span of $\hat U^{(\ell)}$.
Applying this projections to each mode of $\mathcal{T}$ yields the Tucker approximation
\begin{align}
    \mathcal T & \approx  (\mathcal T \times_1 \Phi_{I_1}^T \times_2 \dots \times_d \Phi_{I_d}^T) \times_1 Q^{(1)} (\Phi_{I_1}^T Q^{(1)})^{-1} \times_2 \dots \times_d Q^{(d)} (\Phi_{I_d}^T Q^{(d)})^{-1} \nonumber \\
    &= \mathcal C \times_1 U_1 \times_2 \cdots \times_d U_d, \label{eq:Step1Tucker}
\end{align}
with factor matrices $U^{(\ell)} = Q^{(\ell)} (\Phi_{I_\ell}^T Q^{(\ell)})^{-1}$ and core tensor $\mathcal C = \mathcal{T}(I_1,\dots,I_d)$.
Note that the subtensor $\mathcal{C}$ of $\mathcal{T}$ will not be explicitly formed; we only need to evaluate some its entries when computing the TT approximation of $\mathcal C$ in Section~\ref{sec:PhaseTT}. 
\revision{We want to remark that the core tensor obtained using the HOSVD is not a subtensor of $\mathcal{T}$. }

In~\cite{Stroessner21}, it is shown how the error of the approximation~\eqref{eq:Step1Tucker} is linked to $\norm{\mathcal{T} - \mathcal{T} \times_\ell Q^{(\ell)}(Q^{(\ell)})^T}_F$ through $\norm{(\Phi_{I_\ell}^T Q^{(\ell)})^{-1}}_2$, where $\norm{\cdot}_2$, $\norm{\cdot}_F$ denote the spectral and Frobenius norms, respectively. 
The term $\norm{\mathcal{T} - \mathcal{T} \times_\ell Q^{(\ell)}(Q^{(\ell)})^T}_F$ is small when all mode-$\ell$ fibers of $\mathcal T$ are well approximated by the span of $\hat U^{(\ell)}$. The term $\norm{(\Phi_{I_\ell}^T Q^{(\ell)})^{-1}}_2$ depends on the choice of the index set $I_\ell$.
We apply discrete empirical interpolation~\cite{Chaturantabut10} as formalized in Algorithm~\ref{alg:DEIM} to $Q^{(\ell)}$ to find a suitable index set $I_\ell$.
The whole process of computing $U^{(\ell)}$ and $\mathcal C$ from $\hat U^{(\ell)}$ is summarized in Algorithm~\ref{alg:Tucker}.
Alternative randomized algorithms to compute Tucker approximations have been proposed in~\cite{Minster20,Saibaba21}.

\begin{algorithm}[!ht]
\caption{Discrete empirical interpolation method (DEIM)}\label{alg:DEIM}
\begin{algorithmic}[1]
\State \textbf{Input:} matrix $M \in \mathbb{R}^{n\times m}$ with orthonormal columns
\State \textbf{Output:} index set $I$ if cardinality $m$
\State $I = \{\mathsf{argmax}\ |M(:,1)|\}$
\For $k = 2,\dots,m$
\State $c = M(I,1:k-1)^{-1} M(I,k)$
\State $r = M(:,k) - M(:,1:k-1)c$
\State $I = I \cup \{\mathsf{argmax}\ |r|\}$
\EndFor
\end{algorithmic}
\end{algorithm} 

\begin{algorithm}[!ht]
\caption{Tucker approximation}\label{alg:Tucker}
\begin{algorithmic}[1]
\State \textbf{Input:} procedure to evaluate entries of $\mathcal T \in \R^{(n_1+1)\times \dots \times (n_d+1)}$, tolerance $\eps$, number of samples $s$
\State \textbf{Output:} Matrices $ U^{\ell} \in \R^{(n_\ell +1)\times r_\ell}$ and a procedure to evaluate entries of $\mathcal C \in \R^{r_1\times \dots \times r_d}$, defining a Tucker approximation~\eqref{eq:Step1Tucker}
\For $\ell = 1,\dots,d$
\State $\hat{J}_\ell \gets$ index set $J$ returned by Algorithm~\ref{alg:ACA} applied to $\mathcal{T}^{\{\ell\}}$ with tolerance $\eps$ and $s$ samples.\label{line:ComputeJhat}
\State $\hat U^{(\ell)} = \mathcal{T}^{\{\ell\}}(:,\hat{J}_\ell)$\label{line:ComputeUhat}
\State Compute economic QR decomposition $\hat U^{(\ell)} = Q^{(\ell)} R^{(\ell)}$.
\State $I_\ell = \text{DEIM}(Q^{(\ell)})$ (see Algorithm~\ref{alg:DEIM})
\State $U^{(\ell)} = Q^{(\ell)}(Q^{(\ell)}(I_\ell,:))^{-1}$
\EndFor
\State $\mathcal{C} = \T(I_1,\dots,I_d)$ 
\end{algorithmic}
\end{algorithm} 

\begin{rmrk} \label{rem:Refinement}
For simplicity, we have assumed that the values of $n_1,\dots,n_d$, determining the size of the expansion~\eqref{eq:ChebyshevInterpolantForm} and of $\mathcal T$ are given as input. 
In practice, these values are usually not provided.
To ensure that the approximation error of the final approximation~\eqref{eq:ExtendedFunctionalTT} is small, we use a heuristic to determine $n_1,\dots,n_d$ adaptively. 
We initially set $n_1 = \dots = n_d = 16$. 
After computing $\hat U^{(\ell)}$ in line~\ref{line:ComputeUhat} of Algorithm~\ref{alg:Tucker}, we apply the Chebfun chopping heuristic proposed by Aurentz and Trefethen~\cite{Aurentz17} to $\hat U^{(\ell)}$.
This heuristic decides based on the decay of the coefficients in the columns of the matrix $F^{(\ell)}\hat U^{(\ell)}$ and the tolerance $\eps$ whether $n_\ell$ is sufficiently large. 
If the chopping heuristic indicates that $n_\ell$ is not sufficiently large then one sets $n_\ell \gets 2 n_\ell +1$, updates $\mathcal T$ and repeats the procedure from line~\ref{line:ComputeJhat} to compute a new index set $\hat J_\ell$.
\end{rmrk}

\subsection{TT cores}\label{sec:PhaseTT}

It remains to compute a TT approximation of the core tensor~\eqref{eq:TTofCore} to obtain the desired extended TT approximation~\eqref{eq:ExtendedTTformat}.
For this purpose, we use the rank-adaptive greedy restricted cross interpolation algorithm~\cite{Savostyanov14} as implemented in the TT toolbox\footnote{The TT toolbox is available from \url{https://github.com/oseledets/TT-Toolbox}.}.
This algorithm requires the evaluation of only $\mathcal{O}(drR^2)$ entries of $\mathcal C$  where $R = \max R_\ell$, $r = \max r_\ell$.

The fundamental idea of TT-cross algorithms~\cite{Oseledets10} is to generalize the concept of cross approximation algorithms for matrices, such as Algorithm~\ref{alg:ACA}, to tensors.
\revision{I}ndex sets  $\mathcal I^{\leq \ell} \subset \{ (i_1,\dots,i_{\ell}) | i_k = 1,\dots,r_k,\ k = 1,\dots,\ell\}$  and $\mathcal I^{> \ell} \subset \{ (i_{\ell + 1},\dots,i_{d}) |i_k = 1,\dots,r_k,\ k = \ell+1,\dots,d \}$ each containing $R_{\ell}$ elements for $\ell = 1,\dots,d-1$ \revision{are called nested when $(i_1,\dots,i_{\ell-1},i_{\ell}) \in \mathcal{I}^{\leq \ell} \Rightarrow (i_1,\dots,i_{\ell-1}) \in \mathcal I^{\leq \ell-1}$ and $ (i_{\ell + 1}, i_{\ell+2},\dots,i_{d}) \in \mathcal{I}^{> \ell} \Rightarrow (i_{\ell + 2},\dots,i_{d}) \in \mathcal{I}^{> \ell+1}$. We use a subindex to denote the elements of an index set (sorted in an abritrary but consistent order).
Given nested index sets}, such a \revision{tensor cross} approximation takes the form~\cite{Savostyanov14}
\begin{equation}
\label{eq:TTcrossApprox}
\hat{\mathcal{C}}_{i_1,\dots,i_d}  =
\sum_{s_1 = 1}^{R_1} \sum_{t_1 = 1} ^{R_1} \cdots \sum_{s_{d-1}=1}^{R_{d-1}}
\sum_{t_{d-1}=1}^{R_{d-1}} \prod_{\ell = 1}^d 
\mathcal{C}(\mathcal I^{\leq \ell-1}_{s_{\ell-1}} , i_\ell, \mathcal I^{> \ell}_{t_\ell} )
(C( \mathcal I^{\leq \ell}, \mathcal I^{> \ell} ))^{-1}_{\revision{s_\ell,t_\ell}},
\end{equation}
where we use \rerevision{$\mathcal{I}^{\leq 0} = \mathcal{I}^{> d} = \revision{\emptyset}$}.
Note that this corresponds to a representation in TT format~\eqref{eq:TTformat}, whose TT cores can be computed by evaluating $\mathcal{O}(drR^2)$ entries of $\mathcal{C}$ where $R = \max R_\ell$, $r = \max r_\ell$.
\revision{The error of an approximation of the form~(\ref{eq:TTcrossApprox}), relative to the best approximation error, has been analyzed for nested index sets in~\cite{Savostyanov14} and \emph{without} the restriction to nested index sets in~\cite{Qin2022,Osinsky2019}. In particular,~\cite[Theorem 2]{Osinsky2019} shows that there exists a choice of index sets such that the error in the maximum norm is only a factor $O(R^{1+\lceil \log_2 d \rceil})$ larger than the best approximation error attained by any tensor having the same TT ranks.}

\revision{A practical algorithm is obtained by updating the} index sets sequentially for  $\ell = 1,\dots,d-1$ until we obtain a good approximation. 
For this purpose so-called DMRG supercores~\cite{Savostyanov11} are formed. 
These are defined as subtensors $\mathcal{C}( \mathcal I^{\leq \ell-1}, :,:,\mathcal I^{> \ell +1})$ which are reshaped into tensors in $\R^{R_{\ell-1} \times r_{\ell} \times r_{\ell+1} \times R_{\ell+1}}$. 
The rank-adaptive algorithm proposed by Savostyanov~\cite[Algorithm~2]{Savostyanov14} computes cross approximations of supercores reshaped into matrices \revision{in $\R^{(R_{\ell-1} r_{\ell}) \times (r_{\ell+1}  R_{\ell+1})}$. After reshaping back into tensors, this leads to an approximation} of the form
\[
\mathcal{C}( \mathcal I^{\leq \ell-1} ,:,:, \mathcal{ I}^{> \ell +1}) \approx 
\sum_{s_\ell = 1}^{R_\ell} \sum_{t_\ell = 1}^{R_\ell} 
\mathcal{C}( \mathcal I^{\leq \ell-1} ,:, \mathcal{ J}^{> \ell }_{t_\ell})
(\mathcal{C}( \mathcal J^{\leq \ell}, \mathcal{ J}^{> \ell}))^{-1}_{\revision{s_\ell,t_\ell}}
\mathcal{C}( \mathcal J^{\leq \ell}_{s_\ell}, :, \mathcal{ I}^{> \ell +1}),
\]
where the index sets $\mathcal J^{\leq \ell}, \mathcal J^{\geq \ell}$ are constructed such that $\mathcal I^{\leq \ell} \subset \mathcal J^{\leq \ell}$ and $\mathcal I^{> \ell} \subset \mathcal J^{> \ell}$.
This is achieved by initializing Algorithm~\ref{alg:ACA} with the index sets  $\mathcal I^{\leq \ell},\mathcal I^{> \ell}$ instead of empty sets.
The sampling of the index tuples in Line~\ref{line:tupleSampling} is modified to ensure the indices $\mathcal I^{\leq \ell-1}$ and $\mathcal I^{> \ell+1}$ are nested in $\mathcal J^{\leq \ell}$ and $\mathcal{J}^{> \ell}$ respectively.
Moreover, Algorithm~\ref{alg:ACA}, is stopped after the first iteration, i.e. the index sets
$\mathcal J^{\leq \ell},\mathcal{J}^{> \ell}$ contain at most one element more than $\mathcal I^{\leq \ell},\mathcal{I}^{> \ell}$.
We can now enhance the approximation~\eqref{eq:TTcrossApprox} by replacing the index sets $\mathcal I^{\leq \ell},\mathcal{I}^{> \ell}$ by $\mathcal J^{\leq \ell},\mathcal{J}^{> \ell}$.
Note that this step increases $R_{\ell}$ when $\mathcal I^{\leq \ell} \neq \mathcal J^{\leq \ell}$.
The index update is repeatedly performed for $\ell = 1,\dots, d-1$ and stopped once the error of the approximation~\eqref{eq:TTcrossApprox} is sufficiently small at sample points.
We formalize this procedure in Algorithm~\ref{alg:TTcross}.

\begin{algorithm}[!ht]
\caption{Greedy restricted cross interpolation algorithm}\label{alg:TTcross}
\begin{algorithmic}[1]
\State \textbf{Input:} procedure to evaluate entries of the tensor $\mathcal C \in \R^{r_1\times \dots \times r_d}$, tolerance $\eps$
\State \textbf{Output:} TT approximation~\eqref{eq:TTofCore} of $\mathcal{C}$ with core tensors $\mathcal H^{(\ell)}$
\State Initialize $\mathcal{I}^{\leq \ell}, \mathcal{I}^{>\ell}$ randomly  with only one element each for $\ell = 1,\dots,d$, i.e. $R_1=\cdots=R_{d-1}=1$.
\While \text{error of the approximation~\eqref{eq:TTcrossApprox} of $\mathcal{C}$ is larger than $\eps$ at sample points}
\For $\ell=1,\cdots,d-1$
\State Apply Algorithm~\ref{alg:ACA} to $\mathcal{C}( \mathcal I^{\leq \ell-1},:,:,\mathcal{ I}^{> \ell +1})$ \revision{reshaped into a matrix of size $(R_{\ell-1} r_{\ell}) \times (r_{\ell+1}  R_{\ell+1})$} and update $\mathcal I^{\leq \ell},\mathcal{I}^{> \ell}$ \revision{and} $R_\ell$ as described in Section~\ref{sec:PhaseTT}.
\EndFor
\EndWhile
\State Compute TT cores $\mathcal H^{(\ell)}$ of the TT approximation~\eqref{eq:TTofCore} corresponding to the approximation~\eqref{eq:TTcrossApprox}.
\end{algorithmic}
\end{algorithm} 

\subsection{EFTT approximation algorithm}
We formalize the overall procedure for computing approximations in the EFTT format~\eqref{eq:ExtendedFunctionalTT} in Algorithm~\ref{alg:Approximation}.
Note that $r_1,\dots,r_d$ and $R_1,\dots,R_{d-1}$ are determined adaptively in lines~\ref{line:MultilinearRanks} and~\ref{line:TTranks} \revision{of Algorithm~\ref{alg:Approximation}}.
Let  $R = \max R_\ell$, $r = \max r_\ell$, $n= \max n_\ell$. 
The total number of evaluations of $f$ in Algorithm~\ref{alg:Approximation} is $\mathcal{O}(dr^3+dnr+dsr^2+drR^2)$.
Applying Algorithm~\ref{alg:TTcross} directly to $\mathcal{T}$ requires $\mathcal{O}(dnR^2)$ evaluations and yields an approximation in FTT format~\eqref{eq:TTformat}.

\begin{algorithm}[!ht]
\caption{EFTT approximation}\label{alg:Approximation}
\begin{algorithmic}[1]
\State \textbf{Input:} function $f:[-1,1]^d \to \mathbb{R}$, tolerance $\eps$, number of samples $s$ 
\State \textbf{Output:} TT cores $\mathcal{H}^{(\ell)}$ and procedures for evaluating univariate functions $u_j^{(\ell)}$ defining an approximation of $f$ in the EFTT format~\eqref{eq:ExtendedFunctionalTT}
\State $n_1 = \dots = n_d = 16$ \label{line:setN}
\State Define a procedure to evaluate entries of the evaluation tensor $\mathcal{T}$ defined in~\eqref{eq:EvaluationTensor}. 
\State Apply Algorithm~\ref{alg:Tucker} with tolerance $\eps$ and $s$ samples to $\mathcal T$ to determine $\hat{U}^{(\ell)}$. Simultaneously, we update $n_1,\dots,n_d$ as well as $\mathcal T$ as described in Remark~\ref{rem:Refinement}.  \label{line:MultilinearRanks}
\State Compute the tensors $\mathcal{H}^{(\ell)}$ by applying Algorithm~\ref{alg:TTcross} with tolerance $\eps$ to $\mathcal{C}$. \label{line:TTranks}
\State Create procedures to evaluate univariate functions $u^{(\ell)}_j(x) = \sum_{i = 0}^{n_\ell} \sum_{k=0}^{n_\ell} F_{k,i}^{(\ell)}  U_{i,j}^{(\ell)} T_k(x)$. \label{line:final}
\end{algorithmic}
\end{algorithm} 

\revision{Note that Algorithm~\ref{alg:Approximation} may also be of interest in other situations, not necessarily related to functions. Lines~\ref{line:MultilinearRanks}--\ref{line:final} are of relevance whenever the entries can be (cheaply) evaluated but explicit storage of the tensor as a whole is too expensive.}

%% file: 4NumericalComparison.tex
\section{Numerical experiments} \label{sec:NumericalExperiments}

In this section, we present numerical experiments\footnote{The MATLAB and Python code to reproduce these results is available from \url{https://github.com/cstroessner/EFTT}.} to study the performance of Algorithm~\ref{alg:Approximation}. Unless mentioned otherwise the sample size in Algorithm \ref{alg:ACA} is set to $s = \min\{\bar{n}/2,50\}$, where $\bar{n} = ((n_1+1)(n_2+1)\cdots(n_d+1))^{1/d}$, and all tolerances are set to $10^{-10}$. The approximations error is measured via a Monte-Carlo estimation of the relative $L^2$-error using $10\,000$ samples.

\subsection{Comparison to a direct TT approximation}\label{sec:NumericalDirectTTComparison}
We first compare Algorithm~\ref{alg:Approximation} yielding an approximation in the EFTT format~\eqref{eq:ExtendedFunctionalTT} to a direct TT approximation  of the evaluation tensor $\mathcal{T} \in \R^{100 \times \dots \times 100}$~\eqref{eq:EvaluationTensor} using Algorithm~\ref{alg:TTcross} as in~\cite{Bigoni16}.
The latter approach yields an approximation in the FTT format~\eqref{eq:FunTTformat}.
Since a direct TT approximation is not adaptive with respect to the polynomial degree, we fix the polynomial degree throughout this comparison, i.e., we replace line~\ref{line:setN} by $n_1=\dots=n_d = 100$ and do not update adaptively in line~\ref{line:MultilinearRanks} of Algorithm~\ref{alg:Approximation}.

\paragraph{Benchmark functions.}
For the set of benchmark functions defined in Appendix~\ref{sec:AppendixTestFunctions}, we on average reduce the number of function evaluations required by $30.6\%$ and the required storage by $41.6\%$ when using Algorithm~\ref{alg:Approximation} compared to the direct TT approximation of the coefficient tensor.
For the Ackley function the reduction is $88.8\%$ in terms of function evaluations and $93\%$ in terms of storage.
For the Borehole function, we have $r_\ell = R_\ell^2$, i.e. we can not compress the FTT approximation further.
In this case, the direct TT approximation is slightly more efficient. 
At the same time, this demonstrates that all other test functions, taken from a wide range of applications, can be stored more efficiently using our approach.
Details\revision{, including insights into the reliability of the algorithms,} can be found in Table~\ref{tab:DirectVsEFTT} in the appendix.

\begin{figure}[!ht]
    \centering
    \includegraphics[width=\textwidth]{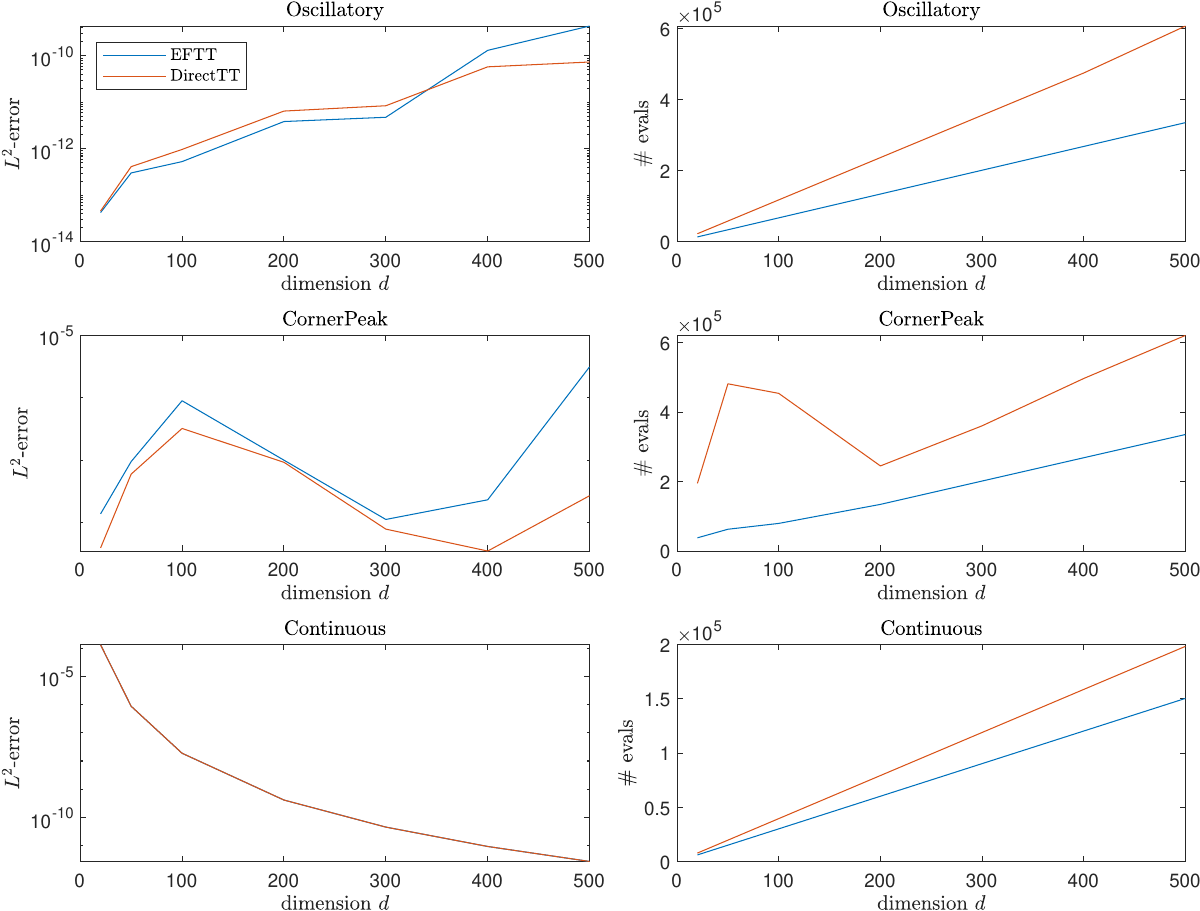}
    \caption{We apply Algorithm~\ref{alg:Approximation} (EFTT) and an approximation of the evaluation tensor using Algorithm~\ref{alg:TTcross} (DirectTT) to approximate the Genz functions (see Appendix~\ref{appendix:Genz}) for dimension $d \in \{20,50,100,200,300,400,500\}$. \revision{Each row corresponds to one test function.} Left: \revision{We plot the $L^2$ error for EFTT and DirectTT}. Right: We plot the \revision{number of} function evaluations required to compute the approximation for EFTT \revision{and} DirectTT. Throughout this figure, we sample $30$ different parameters for the Genz functions as in~\cite{Bigoni16} and display the geometric mean of the error and the arithmetic mean for the number of function evaluations.}
    \label{fig:RankImpact}
\end{figure}

\paragraph{Genz functions}
In order to quantify the impact of the dimension on the approximation error, we study the Genz~\cite{genz1987package} functions defined in Appendix~\ref{appendix:Genz}.
Our numerical results, displayed in Figure~\ref{fig:RankImpact}, demonstrate that Algorithm~\ref{alg:Approximation} requires fewer function evaluation compared to a direct TT-cross approximation. 
However, in very large dimensions ($d>300$) the direct TT-cross approach leads to slightly more accurate approximations. 
\revision{A possible explanation for this observation might be} that, for Algorithm~\ref{alg:Approximation}, the error of $d$ projections needs to be taken into account in addition to the TT-cross approximation error. 
\revision{Since the errors introduced by each projection might affect the overall approximation error multiplicatively (see~\cite[Theorem~2]{Stroessner21}), this could lead to an approximation error growing exponentially in very large dimensions.}
Moreover, approximating the smaller tensor $\mathcal{C}$ instead of $\mathcal{T}$ directly, might lead to slightly worse conditioned inverse matrices in~\eqref{eq:TTcrossApprox}.

\subsection{Comparison to the FTT approximation algorithm}\label{sec:ComapareC3py}
In this section, we compare the performance of our proposed Algorithm \ref{alg:Approximation} based on the EFTT format~\eqref{eq:ExtendedFunctionalTT} to the approximation algorithm proposed in \cite{Gorodetsky19,Gorodetsky17} and implemented in the c3py package \footnote{The c3py package is available from \url{https://github.com/goroda/Compressed-Continuous-Computation}}. 
The c3py algorithm uses a continuous variant of the TT-cross algorithm to compute approximations in the FTT format~\eqref{eq:FunTTformat}.

The c3py package is based on the Legendre polynomials $P_0(x) =1,\ P_1(x) = x,\ (k+2) P_{k+2}(x) = (2k+3) x P_{k+1}(x)-(k+1)P_{k}(x)$ for $k=0,\dots,m-2$~\cite{Koepf98} with degree at most $m \in \mathbb N$ instead of Chebyshev polynomials.
This leads to multivariate polynomial approximations with degree at most $(m_1,\dots,m_d)$ of the form
\begin{equation}
\label{eq:LegendreBasis}
    \sum_{i_1=0}^{m_1} \dots \sum_{i_d=0}^{m_d} \mathcal{B}_{i_1,\dots,i_d} P_{i_1}(x_1)  \cdots P_{i_d}(x_d),
\end{equation}
where $\mathcal{B} \in \R^{(m_1+1) \times \dots \times (m_d+1)}$. 
\revision{The package essentially computes $\mathcal{B}$ in TT format~\eqref{eq:TTformat}.}

\revision{The c3py package obtains the $m$ coefficients for approximating a univariate function using Clenshaw-Curtis quadrature~\cite{Clenshaw60} in $n>m$ points to approximately compute $L^2$-projections onto the Legendre basis functions~\cite{Gorodetsky17}. 
To obtain a fair comparison~\cite{Boyd14}, we modify our approximation algorithm slightly to also include such a Clenshaw-Curtis quadrature. 
Given the evaluation tensor $\mathcal{T}\in \R^{(n_1+1)\times\cdots\times (n_d+1)}$~\eqref{eq:EvaluationTensor}, we compute} 
$\mathcal{B} = \mathcal{T} \times_1 E^{(1)} \times_2 \dots \times_d E^{(d)},$
where
$E^{(\ell)} \in \R^{(n_\ell+1)\times (m_\ell+1)}$ is defined as
\begin{equation*}
    E^{(\ell)} = 
\begin{pmatrix}
\frac{1}{2}w_0^{(\ell)} P_0(x_0^{({ \ell })}) & \frac{1}{2}w_1^{(\ell)} P_0(x_1^{({ \ell })}) & \frac{1}{2}w_2^{(\ell)} P_0(x_2^{({ \ell })}) & \dots & \frac{1}{2}w_{m_{\ell}}^{(\ell)} P_{0}(x_{m_{{ \ell }} }^{({ \ell })}) \\
\frac{3}{2}w_0^{(\ell)} P_1(x_0^{({ \ell })}) & \frac{3}{2}w_1^{(\ell)} P_1(x_1^{({ \ell })}) &  \frac{3}{2}w_2^{(\ell)} P_1(x_2^{({ \ell })}) & \dots & \frac{3}{2}w_{m_{\ell}}^{(\ell)} P_{1}(x_{m_{{ \ell }} }^{({ \ell })}) \\
\frac{5}{2}w_0^{(\ell)} P_2(x_0^{({ \ell })}) &  \frac{5}{2}w_1^{(\ell)} P_2(x_1^{({ \ell })}) &  \frac{5}{2}w_2^{(\ell)} P_2(x_2^{({ \ell })}) & \dots & \frac{5}{2}w_{m_{\ell}}^{(\ell)} P_{2}(x_{m_{{ \ell }} }^{({ \ell })})
\\ \vdots & \vdots & \vdots & \ddots & \vdots \\
\frac{2n_{\ell}+1}{2}w_0^{(\ell)} P_{n_{{ \ell }}}(x_0^{({ \ell })}) & \frac{2n_{\ell}+1}{2}w_1^{(\ell)} P_{n_{{ \ell }}}(x_1^{({ \ell })}) & \frac{2n_{\ell}+1}{2}w_2^{(\ell)} P_{m_{{ \ell }}}(x_2^{({ \ell })}) & \dots & \frac{2n_{\ell}+1}{2}w_{m_{\ell}}^{(\ell)} P_{n_{{ \ell }}}(x_{m_{{ \ell }} }^{({ \ell })})
\end{pmatrix},
\end{equation*}
where the nonnegative weights $w_{i_\ell}^{(\ell)}$ are defined for ${i_\ell} = 0,\ldots,m_{\ell}$, $\ell = 1,\dots,d$ as in \cite{Waldvogel06}:
\begin{align*}
w_{i_\ell}^{(\ell )} &= \frac{c_\ell^{(\ell )}}{(m_\ell )} \big(1 -\sum_{j=1}^{\lfloor m_\ell /2 \rfloor} \frac{b_{j}^{(\ell )}}{4j^2-1}\cos \big(\frac{2j i_\ell \pi}{m_\ell}\big)\big), \nonumber \\
c_{i_\ell}^{(\ell )} &= \left\{
\begin{array}{ll}
1 & {i_\ell}=0\textrm{ or } {i_\ell}=m_{\ell }, \\
2 & \, \textrm{otherwise,} \\
\end{array}
\right. \nonumber \\
b_{i_\ell}^{(\ell )} &= \left\{
\begin{array}{ll}
1 & i_\ell=\frac{m_\ell }{2}, \\
2 & \, \textrm{otherwise}. \\
\end{array}
\right. \nonumber  
\end{align*}
\revision{Note that} the matrix $E^{(\ell)}$ encodes Clenshaw-Curtis quadrature in $n_\ell+1$ nodes to (approximately) compute the Legendre coefficients via $L^2$-projections \revision{as in c3py}.
Using analogous constructions as in Section~\ref{sec:problemsetting}, we can transform an extended TT approximation of the evaluation tensor~\eqref{eq:ExtendedTTformat} into a EFTT approximation with univariate functions represented in terms of linear combinations of Legendre polynomials. \revision{This can also be seen as storing $\mathcal{B}$ in~\eqref{eq:LegendreBasis} in the extended TT format~\eqref{eq:ExtendedTTformat}.}

In the following experiments we set $n_{\ell} = 2m_{\ell}$ for $\ell = 1,\dots,d$ to ensure accurate quadrature. To determine the polynomial degrees $m_\ell$ adaptively (see Remark~\ref{rem:Refinement}), we follow the fiber adaptation strategy of c3py (see \cite[Section 3.6.1]{Gorodetsky17}): We progressively increase the degrees until four sequential Legendre coefficients are smaller than the tolerance $10^{-10}$ or until the maximum of $m_\ell = 105$ has been reached. In the c3py algorithm, we set all tolerances to $10^{-10}$.

\paragraph{Benchmark functions.}
We apply both the c3py algorithm from~\cite{Gorodetsky19} and our novel Algorithm~\ref{alg:Approximation} to approximate the benchmark functions defined in Appendix~\ref{sec:AppendixTestFunctions}. Our numerical results displayed in Table~\ref{tab:C3PYEFTT} in the appendix demonstrate that our novel algorithm typically requires fewer function evaluations and still achieves the same accuracy compared to c3py. For the Ackley function, our approach reduces the number of required function evaluations by more than $96\%$ and the storage by more than $90\%$. At the same time, our approach is slightly more accurate.

\begin{figure}[!ht]
    \centering
    \subfloat{{\includegraphics[width=0.4\textwidth]{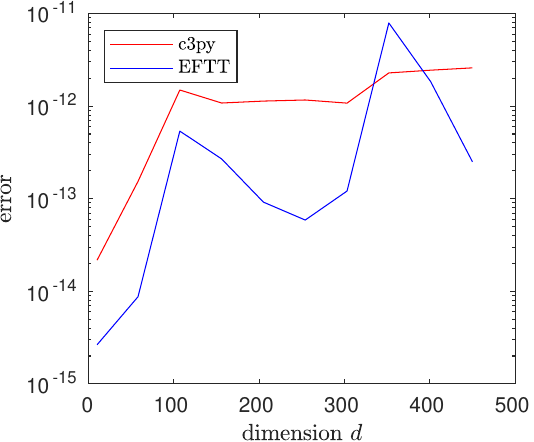}}}%
    \qquad
    \subfloat{{\includegraphics[width=0.4\textwidth]{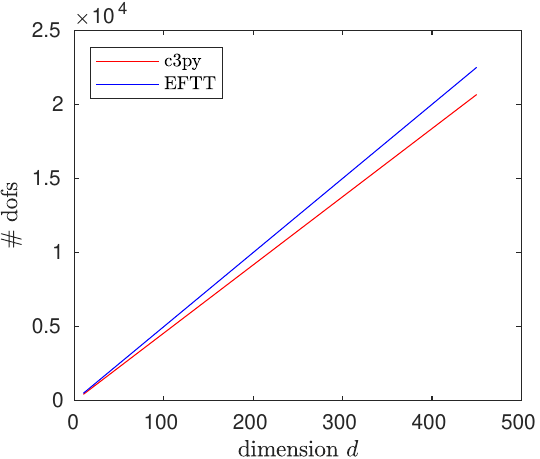}}}%
    \caption{We apply Algorithm~\ref{alg:Approximation} (EFTT) and the algorithm in the c3py package to approximate $f(x_1,\dots,x_d) = \sin(x_1 + x_2 + \dots + x_d)$ by functional low-rank approximation with Legendre polynomial basis functions~\eqref{eq:LegendreBasis} as in~\cite[Figure~3-6]{Gorodetsky17}. Left: We plot the relative error of the integral of the approximations. Right: We plot the number of function evaluations required to compute the approximation. }
    \label{fig:SinIntegration}
\end{figure}

\paragraph{Integration of sin function}
In the following, we repeat the experiment from~\cite[Figure~3-6]{Gorodetsky17}. 
The function $f(x_1,\dots,x_d) = \sin(x_1 + x_2 + \dots + x_d)$ can be represented in FTT format with TT ranks $(2,2,\dots,2)$ and its integral over the domain $[0,1]^d$ is known analytically~\cite{Gorodetsky17}. 
In Figure~\ref{fig:SinIntegration}, we plot the error of the integral of the approximations computed via Algorithm~\ref{alg:Approximation} and the c3py algorithm. 
The figure shows, that our novel approach achieves a similar level of accuracy.
The difference in the number of function evaluations is rather small, since the multilinear ranks and TT ranks of the function are small.

\paragraph{Application: Uncertainty quantification.}
A classical application of multivariate function approximation is the computation of surrogates for uncertainty quantification~\cite{Xiu17,Sudret17}.
For the approximation of the quantity of interest mapping defined in Appendix~\ref{appendix:PDEdefinition}, we find in Table~\ref{tab:Timescompare} that Algorithm~\ref{alg:Approximation} leads to approximations requiring less storage compared to c3py.

\begin{table}[!ht]
\centering
\begin{tabular}{llrrr}
                        &      &  \multicolumn{1}{c}{error}               &  \multicolumn{1}{c}{\# evals} &  \multicolumn{1}{c}{\# dofs} \\ \hline
\multirow{2}{*}{$d=4$}  & c3py & $4.21 \cdot 10^{-5}$ & $360$   & $168$         \\
                        & EFTT & $1.05 \cdot 10^{-3}$  & $337$  & $60$                \\ \hline
\multirow{2}{*}{$d=9$}  & c3py & $9.79\cdot10^{-5}$   & $960$   & $448$         \\
                        & EFTT & $1.30 \cdot 10^{-3}$ & $757$ & $153$               \\ \hline
\multirow{2}{*}{$d=16$} & c3py & $1.73\cdot 10^{-4}$  & $1800$   & $840$         \\
                        & EFTT & $7.62 \cdot 10^{-4}$ & $1345$ & $240$                 \\ \hline
\end{tabular}
\caption{We apply Algorithm~\ref{alg:Approximation} (EFTT) and the algorithm in the c3py package to approximate the quantity of interest map  $Q:[-1,1]^d \to \R$ defined in Appendix~\ref{appendix:PDEdefinition} for $d\in \{4,9,16\}$. All tolerances are set to $10^{-3}$. The table displays the $L^2$-error, the required number of evaluations (evals) and the number of degrees of freedom (dofs) for both approximations.}
\label{tab:Timescompare}
\end{table}

%% file: 5Conclusion.tex
\section{Conclusion and outlook}
For a broad range of functions arising in various applications, including mechanical engineering and uncertainty quantification, our novel algorithm for computing approximations in the EFTT format outperforms existing algorithms based on the FTT format in terms of function evaluations.

In this paper, we have not discussed how to compute numerically with functions once they are compressed in EFTT format. For arithmetic operations, one would need to extend the ideas from~\cite{Trefethen15,Townsend14} to EFTT.  When a function is defined implicitly as the solution of a partial differential equation, extensions of~\cite{Townsend15,CS21preprint} based on ideas from~\cite{Shi21,Kressner09,Psenka20} could be considered.

%% file: 6Apendix.tex
\appendix
\section{Test functions} \label{sec:AppendixTestFunctions}

In Table~\ref{tab:testfuncs1} and Table~\ref{tab:testfuncs2} , we define the test functions for our numerical experiments. 
Note that the functions are defined on different tensor product domains. 
In our experiments, we map the domain of these functions onto $[-1,1]^d$ using an affine linear transformation.

\begin{table}[!ht]
\centering
\begin{tabular}{lccc}
\toprule
\multicolumn{1}{c}{Function}                                                                                                                                                                                                                                   & $d$                 & Domain                        & References \\
\midrule
$\displaystyle \begin{aligned}f_{\text{Ackley}}(\mathbf{x}) =& -20 \exp\big(-0.2 \sqrt{\frac{1}{7}\sum_{i=1}^7 (x_i)^2}\big) \\&- \exp\big(\frac{1}{7} \sum_{i=1}^7 \cos(2 \pi (x_i))\big) + 20 + e^1\end{aligned}$                                        & $7$  & $[-32.768,32.768]^7$    &   \cite{simulationlib,jamil2013literature}   \\
$\displaystyle f_{\text{Alpine}}(\mathbf{x}) = \sum_{i=1}^7 | x_i \sin(x_i) + 0.1 x_1|$                                                                                                                                                                        &  $7$                   & $[-10,10]^7$         &  \cite{jamil2013literature,rahnamayan2007novel} \\
$\displaystyle f_{\text{Dixon}}(\mathbf{x}) = (x_1 -1)^2 + \sum_{i=2}^7 i \cdot (2x_i^2 -x_{i-1})^2$                                                                                                                                                           &                $7$      & $[-10,10]^7$                  &  \cite{jamil2013literature,simulationlib} \\
$\displaystyle f_{\text{Exponential}}(\mathbf{x}) = - \exp\big( -\frac{1}{2} \sum_{i=1}^7 x_i^2 \big)$                                                                                                                                                         &                   $7$   & $[-1,1]^7$        &     \cite{jamil2013literature,rahnamayan2007opposition} \\
$\displaystyle f_{\text{Griewank}}(\mathbf{x}) = \sum_{i=1}^7 \frac{x_i^2}{4000} - \prod_{i=1}^d \cos\big(\frac{x_i}{\sqrt{i}} \big) + 1$                                                                                                                     &                 $7$     & $[-600,600]^7$     &    \cite{simulationlib,jamil2013literature}  \\
$\displaystyle f_{\text{Michalewicz}}(\mathbf{x}) = - \sum_{i=1}^7 \sin\big(x_i\big) \sin^{20}\big( \frac{ix_i^2}{\pi}\big)$                                                                                                                                   &                   $7$   & $[0,\pi]^7$       &     \cite{simulationlib,vanaret2020certified} \\
$\displaystyle f_{\text{Piston}}(\mathbf{x})$                                                                                                                                                                                                                  &                   $7$   & see text                  &  \cite{simulationlib,zankin2018gradient} \\
$\displaystyle f_{\text{Qing}}(\mathbf{x}) = \sum_{i=1}^7 (x_i^2-i)^2$                                                                                                                                                                                         &                 $7$     & $[0,500]^7$                   & \cite{jamil2013literature,qing2005dynamic} \\
$\displaystyle f_{\text{Rastrigin}}(\mathbf{x}) = 70 + \sum_{i=1}^7 (x_i^2 - 10 \cos(2\pi \cdot x_i))$                                                                                                                                                         &                   $7$   & $[-5.12,5.12]^7$              & \cite{dieterich2012empirical,simulationlib} \\
$\displaystyle f_{\text{Rosenbrock}}(\mathbf{x}) = \sum_{i=1}^6 (100 \cdot(x_{i+1} - x_i^2)^2 + (1 - x_i)^2)$                                                                                                                                                  &                    $7$  & $[-2.048,2.048]^7$            &  \cite{jamil2013literature,simulationlib}\\
$\displaystyle f_{\text{Schaffer}}(\mathbf{x}) = \sum_{i=1}^{6} \big(0.5 + \frac{\sin^2\big(\sqrt{x_i^2+x_{i+1}^2}\big)-0.5}{\big(1 + 0.001 (x_i^2+x_{i+1}^2) \big)^2 }\big)$                                                                                  &                  $7$    & $[-100,100]^7$                & \cite{jamil2013literature,simulationlib} \\
$\displaystyle f_{\text{Schwefel}}(\mathbf{x}) = 2932.8803 - \sum_{i=1}^7 x_i \cdot \sin(\sqrt{|x_i|})$                                                                                                                                                        &                  $7$    & $[-500,500]^7$                & \cite{dieterich2012empirical,simulationlib} \\
\bottomrule
\end{tabular}
\caption{Test functions from \cite[Table 1]{Chertkov22}.}
\label{tab:testfuncs1}
\end{table}

\begin{table}[!ht]
\centering
\begin{tabular}{lccc}
\toprule

\multicolumn{1}{c}{Function}                                                                                                                                                                                                                                   & $d$                 & Domain                        & References \\
\midrule

$\displaystyle f_{\text{Borehole}}(\mathbf{x})$                                                                                                                                                                                                                & 8                   & see text &  \cite{simulationlib,an2001quasi} \\
$\displaystyle f_{\text{OTL Circuit}}(\mathbf{x})$                                                                                                                                                                                                             & 6                   &             see text                    &    \cite{simulationlib,moon2012two}  \\
$\displaystyle f_{\text{Robot Arm}}(\mathbf{x})$                                                                                                                                                                                                               & 8                   &      see text                           &  \cite{simulationlib,an2001quasi} \\
$\displaystyle f_{\text{Wing Weight}}(\mathbf{x})$                                                                                                                                                                                                             & 10                  &            see text                     &     \cite{simulationlib,sobester2008engineering} \\
$\displaystyle f_{\text{Friedman}}(\mathbf{x}) = 10 \sin \left(\pi x_{1} x_{2}\right)+20\left(x_{3}-0.5\right)^{2}+10 x_{4}+5 x_{5}$                                                                                                                           & 5                   & $[0,1]^5$                     &     \cite{simulationlib,friedman1991multivariate} \\
$\displaystyle f_{\text{G\&L}}(\mathbf{x}) = \exp \left[\sin \left(\left(0.9\left(x_{1}+0.48\right)\right)^{10}\right)\right]+x_{2} x_{3}+x_{4}$                                                                                                               & 6                   & $[0,1]^6$                     &    \cite{simulationlib,gramacy2009adaptive}  \\
$\displaystyle \begin{aligned}f_{\text{D\&P 8D}}(\mathbf{x}) =& 4\left(x_{1}-2+8 x_{2}-8 x_{2}^{2}\right)^{2}+\left(3-4 x_{2}\right)^{2}\\&+16 \sqrt{x_{3}+1}\left(2 x_{3}-1\right)^{2}+\sum_{i=4}^{8} i \ln \left(1+\sum_{j=3}^{i} x_{j}\right)\end{aligned}$  & 8                   & $[0,1]^8$                     &    \cite{simulationlib,dette2010generalized}  \\
$\displaystyle f_{\text{D\&P Exp}}(\mathbf{x})= 100\left(e^{-2 / x_{1}^{1.75}}+e^{-2 / x_{2}^{1.5}}+e^{-2 / x_{3}^{1.25}}\right)$                                                                                                                              & 3                   & $[0,1]^3$                     &   \cite{simulationlib,dette2010generalized}  \\
\bottomrule
\end{tabular}
\caption{Test functions from \cite[Table 3.2]{Gorodetsky17}.}
\label{tab:testfuncs2}
\end{table}

Some of these functions do not fit into the format of the table. These are defined in the following:
\[
    \displaystyle f_{\text{Piston}}(M,S,V_0,k,P_0,T_a,T_0) = 2\pi
\sqrt{ \frac{M}{k+S^2 \frac{P_0 V_0}{T_0}\frac{T_a}{V^2}}}, \]
where
\[
 V= \frac{S}{2k} \big( \sqrt{A^2 + 4k \frac{P_0V_0}{T_0}T_a - A}   \big) \text{ and }
 A = P_0S + 19.62M- \frac{kV_0}{S},
\]
with $M \in [30,60],\ S \in [0.005,0.02],\ V_0 \in [0.002,0.01],\ k = [1000,5000],\ P_0 \in [90000, 110000],\ T_a \in [290,296],\ T_0 \in [340,360]$.
\[
    \displaystyle f_{\text{Borehole}}(r_w,r,T_u,H_u,T_l,H_l,L,K_w) = 
    \frac{2\pi T_u (H_u-H_l)}{ 
    \log(r/r_w) \big(
    1 + \frac{2LT_u}{\log(r/r_w)r_w^2K_w}+\frac{T_u}{T_l}
    \big),
    }
\]
with $r_w \in [0.05,0.15],\ r \in [100,50000],\ T_u \in [63070,115600],\ H_u \in [990,1110],\ T_l \in [63.1,116],\ H_l \in [700,820],\ L \in [1120,1680],\ K_w \in [9855,12045]$.
\[ 
\displaystyle f_{\text{OTL Circuit}}(b_1,b_2,f,c_1,c_2,\beta) = 
\frac{(\frac{12b_2}{b_1+b_2} + 0.74)\beta(c_2+9)}{\beta(c_2+9)+f}
+ \frac{11.35f}{\beta(c_2+9)+f}+
\frac{0.74 f \beta(c_2 +9)}{(\beta(c_2+9)+f)c_1}
,
\]
with $b_1 \in [50,150],\ b_2 \in [25,70],\ f \in [0.5,3],\ c_1 \in [1.2,2.5],\ c_2 \in [0.25,1.2],\ \beta \in [50,300]$.
\[
    \displaystyle f_{\text{Robot Arm}}(\theta_1,\theta_2,\theta_3,\theta_4,L_1,L_2,L_3,L_4) = 
\sqrt{u^2+v^2},
\]
where $u = \sum_{i=1}^4 L_i \cos(\sum_{j=1}^4 \theta_i)$, $v = \sum_{i=1}^4 L_i \sin(\sum_{j=1}^4 \theta_i)$ and $\theta_i \in [0,2\pi],\ L_i \in [0,1]$ for $i = 1,\dots,4$.
\begin{align*}
\displaystyle & f_{\text{Wing Weight}}(S_w,W_f,A,\Delta,q,\lambda,t_c,N_z,W_d,W_p) = \\ & \quad
0.036 S_w^{0.758}W_f^{0.0035} \big( \frac{A}{\cos^2(\Delta)} \big)^{0.6} q^{0.006} \lambda^{0.04} \big( \frac{100 t_c}{\cos(\Delta)}\big)^{-0.3} (N_zW_d)^{0.49} + S_wW_p,
\end{align*}
with $S_w \in [150,200], W_f \in [220,300],\ A \in [6,10],\ \Delta \in [-10,10],\ q \in [16,45],\ \lambda \in [0.5,1],\ t_c \in [0.08,0.18],\ N_z \in [2.5,6],\ W_d \in [1700,2500],\ W_p \in [0.025,0.08].$

\section{Genz functions} \label{appendix:Genz}
In the following, we define the Genz functions~\cite{genz1987package}, which are frequently used to evaluate function approximation and integration schemes. On the domain $[-1,1]^d$, we consider 
\begin{equation*}
\begin{aligned}
    &f_1(\mathbf{x})=\cos \left(2 \pi w_{1}+\sum_{i=1}^{d} c_{i} \dfrac{x_{i} + 1}{2}\right) & \text{(oscillatory)}\\ 
&f_2(\mathbf{x})=\left(1+\sum_{i=1}^{d} c_{i} \dfrac{x_{i} + 1}{2}\right)^{-(d+1)} & \text{(corner peak)}\\
&f_3(\mathbf{x})=\exp \left(-\sum_{i=1}^{d} c_{i}^{2}\left|\dfrac{x_{i} + 1}{2}-w_{i}\right|\right)& \text{(continuous)}\\
\end{aligned}
\end{equation*}
The parameters $w_i$ and $c_i$ are drawn uniformly from $[0,1]$, where $w_i$ act as a shift for the functions while $c_i$ determines the approximation difficulty of the functions. We normalized $c_i$ such that  
\begin{equation*}
    \sum_{i=1}^{d} |c_{i}| = \dfrac{b}{d^{h}},
\end{equation*}
where the scaling constants $h$ and $b$  are defined for each function in \revision{Table~\ref{tab:GenzConstants}}.
Note that $f_1$ can be represented in FTT format~\eqref{eq:FunTTformat} with $\max R_\ell = 2$.
The function $f_3$ is separable.
For $f_2$, we are not aware of any analytic FTT representation.

\begin{table}[!ht]
\centering
\begin{tabular}{c|rrr}
    & \revision{$f_1$} & \revision{$f_2$} & \revision{$f_3$}  \\ \hline
\revision{$b$} & \revision{284.6} & \revision{185.0} & \revision{2040.0} \\
\revision{$h$} & \revision{1.5}   & \revision{2.0}   & \revision{2.0}   
\end{tabular}
\caption{\revision{Scaling constants for the Genz functions as in~\cite{Bigoni16}.}}
\label{tab:GenzConstants}
\end{table}

\section{Parametric PDE problem}\label{appendix:PDEdefinition}
In the following, we recall the example from~\cite[Section~4]{Kressner11}.
Assume $\sqrt{d} \in \mathbb N$. Let $\Omega = [0,1]^2$. We consider the parametric elliptic PDE
\begin{align}
\label{eq:CookiePDE}
    -\nabla \cdot (a(x,p) \nabla u(x,p)) = 1, & \quad x \in \Omega
\end{align}
with homogeneous Dirichlet boundary conditions and parameter $p \in [-1,1]^{d}$.
We define the piecewise constant coefficient $a(x,p): \Omega \times \R^{d}$ as
\[
a(x,p) = \left\{
\begin{array}{ll}
 1.25 + 0.75 p_{\sqrt{d}(t-1)+s} &  x \in \Omega_{s,t} \\
1 & \, \textrm{otherwise,} \\
\end{array}
\right.
\]
where we denote the disk with radius $\rho = 1/(4\sqrt{d}+2)$ centered around $(\rho(4s-1),\rho(4t-1))$ by $\Omega_{s,t}$ for $s,t = 1,\dots,\sqrt{d}$. 
In our numerical experiments, we approximate the quantity of interest $Q:[-1,1]^{d} \to \R$ defined as
\[
Q(p) = \int_0^1 \int_0^1 u(x,p) \partial x_1 \partial x_2,
\]
where $u(x,p)$ denotes the solution of the PDE~\eqref{eq:CookiePDE} for the given parameter $p \in [-1,1]^{d}$.
For each value of $p$, we solve the resulting PDE using a discretization based on linear finite elements.

\section{Experimental result tables}

\begin{table}[!ht]
\centering
\resizebox{\columnwidth}{!}{%
\begin{tabular}{llrrrrrrrr}
    \multicolumn{1}{c}{{Function}} & \multicolumn{1}{c}{{Algorithm}} & \multicolumn{1}{c}{{Error}} & \multicolumn{1}{c}{{\revision{$\sigma$($\log$(error))}}} & \multicolumn{1}{c}{{\# evals}} & \multicolumn{1}{c}{{\revision{$\sigma$(\# evals)}}} & \multicolumn{1}{c}{{\# dofs}} & \multicolumn{1}{c}{\revision{$\sigma$(\# dofs)}} & \multicolumn{1}{c}{{$\max_\ell R_\ell$}} & \multicolumn{1}{c}{{$\max_\ell r_\ell$}} \\ \hline
    \multirow{2}{*}{Ackley}        & EFTT                            & 1.84e-02                    & 2.06e-05                                           & 63152                          & 3.23e+03                                            & 15949                         & 8.11e+02                                         & 14                                       & 10                                       \\ 
                                   & DirectTT                        & 1.84e-02                    & 1.52e-09                                           & 572531                         & 8.76e+04                                            & 225965                        & 1.34e+04                                         & 18                                       &                                          \\ \hline
    \multirow{2}{*}{Alpine}        & EFTT                            & 5.80e-03                    & 6.98e-15                                           & 4677                           & 9.22e-01                                            & 1448                          & 0.00e+00                                         & 2                                        & 2                                        \\ 
                                   & DirectTT                        & 5.80e-03                    & 2.59e-14                                           & 6860                           & 1.97e+00                                            & 2400                          & 0.00e+00                                         & 2                                        &                                          \\ \hline
    \multirow{2}{*}{Dixon}         & EFTT                            & 1.14e-13                    & 5.27e+00                                           & 11872                          & 2.17e+01                                            & 3548                          & 4.00e+00                                         & 3                                        & 5                                        \\ 
                                   & DirectTT                        & 2.66e-14                    & 1.84e-02                                           & 13022                          & 2.16e+00                                            & 5100                          & 0.00e+00                                         & 3                                        &                                          \\ \hline
    \multirow{2}{*}{Exponential}   & EFTT                            & 2.10e-14                    & 1.32e-03                                           & 2108                           & 0.00e+00                                            & 707                           & 0.00e+00                                         & 1                                        & 1                                        \\ 
                                   & DirectTT                        & 2.09e-14                    & 3.57e-15                                           & 2585                           & 1.62e+00                                            & 700                           & 0.00e+00                                         & 1                                        &                                          \\ \hline
    \multirow{2}{*}{Griewank}      & EFTT                            & 1.92e-07                    & 5.96e-01                                           & 8089                           & 2.68e+01                                            & 2252                          & 4.00e+00                                         & 3                                        & 3                                        \\ 
                                   & DirectTT                        & 1.54e-07                    & 2.23e-10                                           & 13023                          & 2.38e+00                                            & 5100                          & 0.00e+00                                         & 3                                        &                                          \\ \hline
    \multirow{2}{*}{Michalewicz}   & EFTT                            & 4.05e-02                    & 2.80e-15                                           & 4677                           & 8.69e-01                                            & 1448                          & 0.00e+00                                         & 2                                        & 2                                        \\ 
                                   & DirectTT                        & 4.05e-02                    & 1.01e-15                                           & 6860                           & 2.10e+00                                            & 2400                          & 0.00e+00                                         & 2                                        &                                          \\ \hline
    \multirow{2}{*}{Piston}        & EFTT                            & 3.32e-09                    & 2.32e-01                                           & 203484                         & 6.39e+03                                            & 74228                         & 2.13e+03                                         & 24                                       & 11                                       \\ 
                                   & DirectTT                        & 2.93e-09                    & 2.50e-01                                           & 992566                         & 6.91e+03                                            & 412603                        & 3.08e+03                                         & 18                                       &                                          \\ \hline
    \multirow{2}{*}{Qing}          & EFTT                            & 1.09e-13                    & 1.53e+00                                         & 5482                           & 1.39e+00                                            & 2172                          & 0.00e+00                                         & 2                                        & 3                                        \\ 
                                   & DirectTT                        & 2.29e-14                    & 5.70e-03                                           & 6860                           & 2.11e+00                                            & 2400                          & 0.00e+00                                         & 2                                        &                                          \\ \hline
    \multirow{2}{*}{Rastrigin}     & EFTT                            & 2.28e-14                    & 2.70e-03                                           & 4677                           & 8.64e-01                                            & 1448                          & 0.00e+00                                         & 2                                        & 2                                        \\ 
                                   & DirectTT                        & 2.30e-14                    & 5.08e-03                                           & 6860                           & 1.87e+00                                            & 2400                          & 0.00e+00                                         & 2                                        &                                          \\ \hline
    \multirow{2}{*}{Rosenbrock}    & EFTT                            & 2.83e-14                    & 3.71e-01                                         & 10970                          & 1.96e+00                                            & 2798                          & 0.00e+00                                         & 3                                        & 4                                        \\ 
                                   & DirectTT                        & 2.64e-14                    & 1.03e-02                                           & 13023                          & 2.22e+00                                            & 5100                          & 0.00e+00                                         & 3                                        &                                          \\ \hline
    \multirow{2}{*}{Schaffer}      & EFTT                            & 6.75e-02                    & 5.53e-03                                           & 1061290                        & 1.61e+05                                            & 288167                        & 5.11e+04                                         & 39                                       & 40                                       \\ 
                                   & DirectTT                        & 6.73e-02                    & 1.67e-12                                           & 1513169                        & 2.93e+04                                            & 767463                        & 6.73e+03                                         & 30                                       &                                          \\ \hline
    \multirow{2}{*}{Schwefel}      & EFTT                            & 6.58e-04                    & 2.98e-14                                           & 4677                           & 8.73e-01                                            & 1448                          & 0.00e+00                                         & 2                                        & 2                                        \\ 
                                   & DirectTT                        & 6.58e-04                    & 3.91e-14                                           & 6860                           & 2.09e+00                                            & 2400                          & 0.00e+00                                         & 2                                        &                                          \\ \hline
    \multirow{2}{*}{Borehole}      & EFTT                            & 3.95e-02                    & 1.01e-10                                           & 14186                          & 6.84e+02                                            & 3243                          & 9.90e+01                                         & 2                                        & 4                                        \\ 
                                   & DirectTT                        & 3.95e-02                    & 6.47e-11                                           & 10042                          & 5.42e+02                                            & 2318                          & 7.16e+01                                         & 2                                        &                                          \\ \hline
    \multirow{2}{*}{OTL Circuit}   & EFTT                            & 3.71e-11                    & 1.44e+00                                           & 16065                          & 5.21e+02                                            & 3280                          & 7.00e+01                                         & 5                                        & 5                                        \\ 
                                   & DirectTT                        & 8.49e-12                    & 2.51e-01                                           & 27764                          & 3.22e+00                                            & 8300                          & 0.00e+00                                         & 4                                        &                                          \\ \hline
    \multirow{2}{*}{Robot Arm}     & EFTT                            & 7.00e-02                    & 3.31e-01                                           & 500591                         & 1.63e+05                                            & 101847                        & 5.69e+04                                         & 33                                       & 33                                       \\ 
                                   & DirectTT                        & 6.52e-02                    & 5.77e-01                                           & 734573                         & 4.78e+05                                            & 383466                        & 2.52e+05                                         & 34                                       &                                          \\ \hline
    \multirow{2}{*}{Wing Weight}   & EFTT                            & 3.73e-14                    & 2.79e-02                                          & 6692                           & 1.21e+00                                            & 2072                          & 0.00e+00                                         & 2                                        & 2                                        \\ 
                                   & DirectTT                        & 8.29e-14                    & 1.68e-01                                           & 10440                          & 2.22e+00                                            & 3600                          & 0.00e+00                                         & 2                                        &                                          \\ \hline
    \multirow{2}{*}{Friedman}      & EFTT                            & 4.41e-10                    & 4.29e+00                                           & 12317                          & 7.05e+01                                            & 2377                          & 1.80e+01                                         & 4                                        & 4                                        \\ 
                                   & DirectTT                        & 8.84e-12                    & 1.92e+00                                           & 14676                          & 7.92e+02                                            & 3142                          & 1.05e+02                                         & 3                                        &                                          \\ \hline
    \multirow{2}{*}{G \& L}        & EFTT                            & 2.52e-05                    & 1.12e-12                                           & 3278                           & 4.56e-01                                            & 1034                          & 0.00e+00                                         & 2                                        & 2                                        \\ 
                                   & DirectTT                        & 2.52e-05                    & 4.46e-13                                           & 6651                           & 2.27e+00                                            & 1800                          & 0.00e+00                                         & 2                                        &                                          \\ \hline
    \multirow{2}{*}{G \& P 8D}     & EFTT                            & 3.08e-11                    & 4.59e-01                                           & 39724                          & 1.22e+03                                            & 8138                          & 1.96e+02                                         & 7                                        & 7                                        \\ 
                                   & DirectTT                        & 2.64e-11                    & 2.98e-01                                           & 74942                          & 4.79e+03                                            & 30140                         & 1.26e+03                                         & 5                                        &                                          \\ \hline
    \multirow{2}{*}{D \& P Exp}    & EFTT                            & 1.56e-14                    & 3.68e-03                                           & 1990                           & 0.00e+00                                            & 616                           & 0.00e+00                                         & 2                                        & 2                                        \\ 
                                   & DirectTT                        & 1.55e-14                    & 3.43e-04                                           & 2087                           & 1.14e+00                                            & 800                           & 0.00e+00                                         & 2                                        &                                          \\ \hline
\end{tabular}
    }
\caption{As described in Section~\ref{sec:NumericalDirectTTComparison}, we apply Algorithm~\ref{alg:Approximation} (EFTT) and an approximation of the evaluation tensor using Algorithm~\ref{alg:TTcross} (DirectTT) to approximate the test functions defined in Appendix~\ref{sec:AppendixTestFunctions}). For each function, we display the estimated $L^2$ error of the approximations, the number of function evaluations required to construct the approximation, the degrees of freedom in the approximation and the largest $R_\ell$ and $r_\ell$ of the resulting approximation. 
The numbers displayed are the mean (geometric mean for the error and arithmetic mean for all other quantities) after approximating each function $100$ times. \revision{In addition, we display the standard deviation of the logarithm of the error  $\sigma$($\log$(error)), as well as the standard deviation for the number of function evaluations ($\sigma$(\#evals)) and the number of degrees of freedom ($\sigma$(\#dofs)).}}
\label{tab:DirectVsEFTT}
\end{table}

\begin{table}[!ht]
\centering
\begin{tabular}{llrrrrrr}

\multicolumn{1}{c}{Function} & \multicolumn{1}{c}{Algorithm} & \multicolumn{1}{c}{Error} & \multicolumn{1}{c}{\# evals} & \multicolumn{1}{c}{\# dofs} & $\max_\ell m_\ell$ & \multicolumn{1}{c}{$\max_\ell R_\ell$} & \multicolumn{1}{c}{$\max_\ell r_\ell$} \\ \hline
\multirow{2}{*}{Ackley}        & EFTT                        & 1.22e-02                   & 67107                        & 18465                        & 105                & 15                                      & 9                                       \\  
                               & c3py                        & 1.81e-02                   & 2232528                      & 197760                       & 105                & 16                                      &                                         \\ \hline
\multirow{2}{*}{Alpine}        & EFTT                        & 4.08e-03                   & 5464                         & 1518                         & 105                & 2                                       & 2                                       \\  
                               & c3py                        & 6.43e-03                   & 50656                        & 2520                         & 105                & 2                                       &                                         \\ \hline
\multirow{2}{*}{Dixon}         & EFTT                        & 3.21e-14                   & 13756                        & 3714                         & 105                & 3                                       & 5                                       \\  
                               & c3py                        & 2.89e-14                   & 13752                        & 435                          & 21                 & 3                                       &                                         \\ \hline
\multirow{2}{*}{Exponential}   & EFTT                        & 4.82e-15                   & 946                          & 196                          & 27                 & 1                                       & 1                                       \\  
                               & c3py                        & 1.71e-14                   & 10518                        & 133                          & 19                 & 1                                       &                                         \\ \hline
\multirow{2}{*}{Griewank}     & EFTT                        & 8.21e-08                   & 9139                         & 2358                         & 105                & 3                                       & 3                                       \\  
                               & c3py                        & 3.52e-06                   & 51466                        & 4459                         & 105                & 3                                       &                                         \\ \hline
\multirow{2}{*}{Michalewicz}   & EFTT                        & 2.54e-02                   & 5464                         & 1518                         & 105                & 2                                       & 2                                       \\  
                               & c3py                        & 1.45e-01                   & 48745                        & 2443                         & 105                & 2                                       &                                         \\ \hline
\multirow{2}{*}{Piston}        & EFTT                        & 3.71e-09                   & 174188                       & 69019                        & 33                 & 23                                      & 11                                      \\  
                               & c3py                        & 3.85e-05                   & 251760                       & 66080                        & 35                 & 24                                      &                                         \\ \hline
\multirow{2}{*}{Qing}          & EFTT                        & 5.54e-13                   & 6996                         & 2277                         & 105                & 2                                       & 3                                       \\  
                               & c3py                        & 2.86e-14                   & 11776                        & 136                          & 21                 & 2                                       &                                         \\ \hline
\multirow{2}{*}{Rastrigin}     & EFTT                        & 1.91e-14                   & 5463                         & 1518                         & 105                & 2                                       & 2                                       \\  
                               & c3py                        & 1.86e-10                   & 22288                        & 1342                         & 63                 & 2                                       &                                         \\ \hline
\multirow{2}{*}{Rosenbrock}    & EFTT                        & 1.25e-14                   & 4106                         & 835                          & 27                 & 3                                       & 4                                       \\  
                               & c3py                        & 9.43e-14                   & 11530                        & 633                          & 14                 & 3                                       &                                         \\ \hline
\multirow{2}{*}{Schaffer}      & EFTT                        & 7.19e-02                   & 173787                       & 35016                        & 105                & 16                                      & 17                                      \\  
                               & c3py                        & 1.22e-01                   & 3465760                      & 214200                       & 105                & 20                                      &                                         \\ \hline
\multirow{2}{*}{Schwefel}      & EFTT                        & 4.00e-04                   & 5463                         & 1518                         & 105                & 2                                       & 2                                       \\  
                               & c3py                        & 5.45e-04                   & 50656                        & 2496                         & 104                & 2                                       &                                         \\ \hline
\multirow{2}{*}{Borehole}      & EFTT                        & 3.95e-02                   & 6552                         & 1116                         & 32                 & 2                                       & 4                                       \\  
                               & c3py                        & 2.08e-03                   & 14346                        & 577                          & 70                 & 2                                       &                                         \\ \hline
\multirow{2}{*}{OTL Circuit}   & EFTT                        & 7.93e-11                   & 6670                         & 1083                         & 27                 & 5                                       & 5                                       \\  
                               & c3py                        & 4.07e-08                   & 15674                        & 1782                         & 28                 & 5                                       &                                         \\ \hline
\multirow{2}{*}{Robot Arm}     & EFTT                        & 8.12e-02                   & 499954                       & 54760                        & 94                 & 12                                      & 27                                      \\  
                               & c3py                        & 3.85e-01                   & 2018017                      & 228439                       & 105                & 20                                      &                                         \\ \hline
\multirow{2}{*}{Wing Weight}   & EFTT                        & 2.83e-14                   & 2867                         & 560                          & 24                 & 2                                       & 2                                       \\  
                               & c3py                        & 2.15e-13                   & 12224                        & 554                          & 19                 & 2                                       &                                         \\ \hline
\multirow{2}{*}{Friedman}      & EFTT                        & 2.16e-11                   & 5238                         & 404                          & 19                 & 3                                       & 4                                       \\  
                               & c3py                        & 8.08e-05                   & 12142                        & 710                          & 15                 & 4                                       &                                         \\ \hline
\multirow{2}{*}{G \& L}        & EFTT                        & 4.95e-06                   & 1547                         & 356                          & 29                 & 2                                       & 2                                       \\  
                               & c3py                        & 3.51e-02                   & 13928                        & 374                          & 105                & 2                                       &                                         \\ \hline
\multirow{2}{*}{G \& P 8D}     & EFTT                        & 4.77e-11                   & 19527                        & 3902                         & 24                 & 6                                       & 7                                       \\  
                               & c3py                        & 9.54e-10                   & 27336                        & 5136                         & 21                 & 7                                       &                                         \\ \hline
\multirow{2}{*}{D \& P Exp}    & EFTT                        & 1.13e-14                   & 2404                         & 646                          & 105                & 2                                       & 2                                       \\  
                               & c3py                        & 4.78e-10                   & 12162                        & 336                          & 49                 & 2                                       &                                         \\ \cline{1-8} 
\end{tabular}
\caption{As described in Section~\ref{sec:ComapareC3py}, we apply Algorithm~\ref{alg:Approximation} (EFTT) and the algorithm in the c3py package to approximate the test functions defined in Appendix~\ref{sec:AppendixTestFunctions} by functional low-rank approximation with Legendre polynomial basis functions~\eqref{eq:LegendreBasis}. 
For each function, we display the estimated $L^2$-error of the approximations, the number of function evaluations required to construct the approximation, the degrees of freedom in the approximation and the largest $m_\ell$, $R_\ell$ and $r_\ell$ of the resulting approximation.}
\label{tab:C3PYEFTT}
\end{table}